\documentclass[preprint,12pt]{elsarticle}
\usepackage{amsfonts,amssymb,amsmath}
\usepackage[mathscr]{eucal}
\usepackage{geometry,graphicx,enumerate,color}

\begin{document}
\begin{frontmatter}

\title{On typical properties of Hilbert space operators}

\author{Tanja Eisner\fnref{Afn}}
\ead{t.eisner@uva.nl}
\ead[url]{http://staff.science.uva.nl/~eisner/}
\address{KdV Institute for Mathematics, University of Amsterdam, P.O. Box 94248, 1090 GE Amsterdam, The Netherlands.}
\fntext[Afn]{The author was supported by the European Social Fund
and by the Ministry of Science, Research and the Arts
Baden-W\"urttemberg.}

\author{Tam\'as M\'atrai\corref{Tfn}}
\ead{tmatrai@math.utoronto.ca}
\ead[url]{http://www.renyi.hu/~matrait}
\address{University of Toronto, Mathematics Department BA 6290, 40 St. George St., M5S 2E4 Toronto, Ontario,  Canada.}
\cortext[Tfn]{The author was partially supported  by the OTKA
grants K 61600, K 49786 and K 72655, by the NSERC grants 129977,
A-7354, RGPIN 3185-10 and 402762, and by the NSF Grant DMS
0600940.}

\begin{abstract}
We study the typical behavior of bounded linear operators on
infinite dimensional complex separable Hilbert spaces in the norm,
strong-star, strong, weak polynomial and weak topologies. In
particular, we investigate typical spectral properties, the
problem of unitary equivalence of typical operators, and their
embeddability into $C_{0}$-semigroups. Our results provide
information on the applicability of Baire category methods in the
theory of Hilbert space operators.
\end{abstract}


\end{frontmatter}

\def\es{\emptyset}
\def\ve{\varepsilon}
\def\vp{\varphi}
\def\vt{\vartheta}
\def\sm{\setminus}
\def\real{\mathbb{R}}
\def\bO{\mathbb{O}}
\def\bU{\mathbb{U}}
\def\bM{\mathbb{M}}
\def\bW{\mathbb{W}}
\def\bF{\mathbb{F}}
\def\bV{\mathbb{V}}
\def\cmpl{\mathbb{C}}
\def\cl{\textrm{cl}}
\def\bel{\textrm{int}}
\def\egesz{\mathbb{Z}}
\def\nat{\mathbb{N}}
\def\I{\mathbb{I}}
\def\b{\left}
\def\j{\right}
\def\<{\left<}
\def\>{\right>}
\def\rar{\rightarrow}
\def\lrar{\leftrightarrow}
\def\ben{\begin{enumerate}}
\def\een{\end{enumerate}}
\def\Keq{\begin{equation}}
\def\Zeq{\end{equation}}
\def\bit{\begin{itemize}}
\def\eit{\end{itemize}}
\def\pa{\partial}
\def\d{\textrm{d}}
\def\w{\omega}
\def\om{\omega}
\def\W{\Omega}
\def\cpl{\complement}
\def\cf{\textrm{cf}}
\def\card{\textrm{card}}
\def\ss{\subseteq}
\def\sqss{\sqsubseteq}
\def\Dh{\Delta_{h}}
\def\vv{\vert \vert}
\def\mc{\mathcal}
\def\DS{\displaystyle}
\def\TS{\textstyle}
\def\BS{{\bf \Sigma}}
\def\BD{{\bf \Delta}}
\def\BP{{\bf \Pi}}
\def\sign{\textrm{sign}}
\def\Ran{\mathrm{Ran\,}}
\def\gotF{\mathfrak{F}}
\def\gotJ{\mathfrak{J}}
\def\gotI{\mathfrak{I}}
\def\gotU{\mathfrak{U}}
\def\diam{\textrm{diam}}
\def\lin{\mathrm{lin}}
\def\dim{\mathrm{dim~}}
\def\supp{\mathrm{supp~}}
\def\dist{\textrm{dist}}
\def\dom{\mathrm{dom}}
\def\Ker{\mathrm{Ker\,}}
\def\re{\mathrm{Re}}
\def\Id{\mathrm{Id}}
\def\lev{\mathrm{lev}}
\def\bel{\mathrm{int}}
\def\PR{\textrm{Pr}}
\def\bs{\blacksquare}
\def\prf{\textbf{Proof. }}
\def\gotT{\mathfrak{T}}
\def\iff{\Leftrightarrow}

\def\span#1{\mathrm{span}\{ {#1} \}}
\def\sp#1{\langle {#1} \rangle}

\newtheorem{thm}{Theorem}[section]
\newtheorem{defin}[thm]{Definition}
\newtheorem{lem}[thm]{Lemma}
\newtheorem{prob}[thm]{Problem}
\newtheorem{prop}[thm]{Proposition}
\newtheorem{cor}[thm]{Corollary}

\section{Introduction}

Given a property $\Phi$ on the points of a Baire space $X$, we say
that \emph{a typical point of $X$ satisfies $\Phi$}, or simply
that \emph{$\Phi$ is typical}, if the set $\{x \in X\colon x
\text{ satisfies } \Phi\}$ is co-meager in $X$, i.e.\ if $\{x \in
X\colon x \text{ does not satisfy } \Phi\}$ is of first category
in $X$. Many important and classical results in analysis are
concerned with typical properties in particular topological
spaces. Examples include the Banach--Mazurkiewicz theorem (see
e.g.\  \cite{Ba} and \cite{Ma}) stating that the set of continuous
nowhere differentiable functions are residual in
$(C([0,1]),\|\cdot\|_{\infty})$ (see also \cite{B} for a primer on
typical properties of continuous functions),  or the famous result
by P.\ R.\ Halmos \cite{H} and V.\ A.\ Rohlin \cite{Ro}  in
ergodic theory on the existence of weakly mixing but not strongly
mixing transformations.

In this paper we continue the investigations of the first author
and study the typical properties of contractive linear operators
on infinite dimensional complex separable Hilbert spaces in the
norm, strong-star, strong, weak polynomial and weak topologies
(for the definitions,  see Definition \ref{topol}). Typical
properties of various classes of operators have been studied
previously: see e.g.\ \cite{CN} for typical properties of measure
preserving transformations; \cite{BK} and \cite{Ku} for typical
mixing properties of Markov semigroups; \cite{E1}, \cite{ES1} and
\cite{ES3} for typical stability properties; \cite{JS}, \cite{RMS}
and \cite{Si} for typical spectral properties in various very
special families of operators. Our research is different from
these works in several respects. We study typical properties of
contractions as a whole, and we carry out our analysis in several
topologies. Surprisingly, in contrast to classical results, we
mostly obtain ``good'' properties as being typical, and it turns
out that the typical properties may change drastically if the
reference topology is changed.

We obtain the following results. In Section \ref{ssweak}, we
recall some results obtained by the first author (see \cite{E0},
\cite{E2}) about typical properties in the weak topology. In this
topology, a typical contraction is unitary, it has maximal
spectrum and empty point spectrum, it can be embedded into a
$C_{0}$-semigroup, and typical contractions are not unitarily
equivalent. Our results make use of the theory of typical
properties of measures developed by M.\ G.\ Nadkarni
\cite[Chapter 8]{N} (see also \cite{CN}). The importance of the
weak topology in operator theory is an obvious motivation for our
 investigations.

In Section \ref{sspw} we consider the weak polynomial topology.
Our main observations are that the contractions endowed with this
topology form a Polish space, where the set of unitary operators
is a co-meager subset. Since on the set of unitary operators the
weak and weak polynomial topologies coincide, we conclude that the
typical properties of contractions in the weak and weak polynomial
topologies coincide. This part of our work is motivated by the
increasing interest in this unusual topology.

Section \ref{ssstrong} treats the strong topology. We show that a
typical contraction is unitarily equivalent to the infinite
dimensional backward unilateral shift operator. An analogous
result for strongly continuous semigroups is obtained as well. In
particular, the point spectrum of a typical contraction is the
open unit disk, typical contractions are unitarily equivalent and
can be embedded into a $C_{0}$-semigroup. Contrast these results
with the behavior in the weak topology. Just as for the weak
topology, our interest in the strong topology necessitates no
clarification.

We study the strong-star topology in Section \ref{ssstrongstar}.
Our main result gives that the theory of typical properties of
contractions in the strong-star topology can be reduced to the
theories of typical properties of unitary and positive
self-adjoint operators in the strong topology. As a corollary, we
obtain that a typical contraction has maximal spectrum and empty
point spectrum, and two typical contractions are not unitarily
equivalent. We included the strong-star topology in our research
because it may play the most important role while extending our
investigations into general Banach$^{\star}$-algebras.

Section \ref{ssnorm} contains our results on the norm topology,
the only non-separable topology we consider. We obtain that in
this topology there are no such non-trivial typical structural
properties as for the separable topologies. Intuitively, the
reason for this phenomenon is that the norm topology is fine
enough to allow for the coexistence of many different properties
on non-meager sets. We close the paper with an outlook to typical
properties of operators on general Banach spaces, and with a list
of open problems.

Before turning our attention to the proofs, let us justify our
settings. We choose contractions as the underlying set of
operators because it becomes a Baire space in all the five
topologies we consider. By writing the set of bounded linear
operators as a countable union of scaled copies of the set of
contractions, suitable extensions of our result can be easily
obtained. In our investigations of contractions, several other
important classes of operators (e.g.\ isometries, positive
self-adjoint operators, unitary operators) come into play, and we
obtain information about typical properties in these subclasses,
as well.

We work only in infinite dimensional separable Hilbert spaces
because removing any of these assumptions invalidates most of our
results. Infinite dimensionality guarantees that the families of
operators under consideration are sufficiently rich. Separability
is essential for our descriptive set theoretic arguments. The
Hilbert space structure not only allows us to use a well-developed
spectral theory, but also facilitates the construction of
operators using orthogonal decomposition.  It is of limited
importance that our Hilbert spaces are over the complex field;
analogous results  hold in real Hilbert spaces, as well. We defer
the further discussion of possible extensions of our work until
Section \ref{sgenB}.

To conclude this introduction we note that in research works
studying contractions, it is customary to point out that
contractions in general are hard to study.  The theory of
contractions as a whole is often contrasted to the theories of
normal, self-adjoint or unitary operators where a satisfactory
classification can be obtained, e.g.\ via spectral measures (see
e.g.\ \cite{N} and \cite{Y}). Our results provide an explanation
of this intuitive observation. As we pointed out above, as far as
Baire category methods are concerned, in the four separable
topologies we consider, the theory of contractions is reduced to
the theory of unitary operators (weak and weak polynomial
topologies), to the theory of one shift operator (strong topology)
or to the theories of unitary and positive self-adjoint operators
(strong-star topology). Thus if the oversimplified pictures
captured by these separable topologies are dissatisfactory for an
analyst, then necessarily the very fine norm topology has to be
used,  in which case non-separability can be made responsible for
being complicated. Since only such properties can be studied using
Baire category arguments which are non-trivial in a suitable Baire
topology, i.e.\ which hold at least on a non-meager set, our
results outline a limitation of Baire category methods in operator
theory.

\medskip

\textbf{Acknowledgment.} We thank Mariusz Lema\'nczyk and Nicolas
Monod for helpful discussions and for calling our attention to the
problem of unitary equivalence of typical operators, as well as
the referee for careful reading and helpful suggestions.

%\section{Typical properties in various topologies}\label{stypprop}
\section{Preliminaries}\label{stypprop}

As general references, see  \cite{K} for descriptive set theory,
\cite{Y} for functional analysis, and \cite{EN} for semigroup
theory. Recall that a topological space $X$ is a \emph{Baire
space} if every non-empty open set in $X$ is non-meager, or
equivalently if the intersection of countably many dense open sets
in $X$ is dense (see e.g.\ \cite[(8.1) Proposition and (8.2)
Definition p.\ 41]{K}). Polish spaces, i.e.\  separable complete
metric spaces, are well-known examples of Baire spaces.

A set in a topological space is $G_{\delta}$ if it can be obtained
as an intersection of countably many open sets. We will often use
the observation that in every topological space, every dense
$G_{\delta}$ set is co-meager.

Let $X$ be a Baire space and let $\Phi$ be a property on the
points of $X$. We say that \emph{$\Phi$ is a typical property on
$X$}, or that \emph{a typical element of $X$ satisfies $\Phi$} if
$\{x \in X \colon x \textrm{ satisfies }\Phi\}$ is a co-meager
subset of $X$. Note that if $\Phi_{n}$ $(n \in \nat)$ are typical
properties on $X$ then a typical element of $X$ satisfies all
$\Phi_{n}$ $(n \in \nat)$ simultaneously.

In the sequel $(H,\| \cdot \|)$ always denotes an infinite
dimensional complex separable Hilbert space. The scalar product on
$H$ is denoted by $\langle\cdot,\cdot \rangle$. For every $U \ss
H$, $\span{U}$ denotes the linear subspace of $H$ generated by
$U$, and $U^{\bot} = \{x \in H \colon \sp{x,u} = 0~(u \in U)\}$.
For every $U,V \ss H$, we write $U \bot V$ if $\sp{u,v} = 0$ $(u
\in U, v \in V)$.  If $V \leq H$ is a subspace, we define $B_{V} =
\{v \in V \colon \|v\| \leq 1\}$, $S_{V} = \{v \in V \colon \|v\|
= 1\}$.

Let $B(H)$,  $C(H)$, $U(H)$ and $P(H)$ denote the sets of bounded,
contractive, unitary and contractive positive self-adjoint linear
$H \rar H$ operators. The identity operator is denoted by $\Id$.
For every $V \leq H$, the orthogonal projection onto $V$ is
denoted by $\Pr_{V}$. For every $A \in B(H)$ the adjoint of $A$ is
denoted by $A^{\star}$. For every $A \in B(H)$ and $U \ss H$ we
set $A[U]=\{Ax \colon x \in U\}$, $\ker A = \{x \in H \colon Ax =
0\}$ and $\Ran A = A[H]$. For every $A \in B(H)$, let $\mc{O}(A) =
\{UAU^{-1} \colon U \in U(H)\}$.

For every $A \in B(H)$ the spectrum, the point spectrum, the
continuous spectrum, and the residual spectrum of $A$ is denoted
by $\sigma(A)$, $P_{\sigma}(A)$, $C_{\sigma}(A)$, and
$R_{\sigma}(A)$ (see e.g.\ \cite[Definition p.\ 209]{Y}).

We recall the definitions and some elementary properties of the
weak, weak polynomial, strong and strong-star topologies (see
e.g.\ \cite[Definition 3 p.\ 112]{Y}, \cite[Definition p.\
1142]{P}, \cite[Definition p.\ 69]{Y}, and \cite[Definition p.\
220]{W}).

\begin{defin}\label{topol}\rm
Let  $A, A_{n} \in B(H)$ $(n \in \nat)$ be arbitrary.
\ben
\item  We say that \emph{$\{A_{n} \colon n \in \nat\}$ converges to
    $A$ weakly}, $A = w$-$\lim_{n \in \nat} A_{n}$ in notation, if for
  every $x,y\in H$, $\lim_{n \in \nat} \sp{A_{n} x,y} =
\sp{Ax,y}$. The topology corresponding to this notion of
convergence is called the \emph{weak topology}. Topological
notions referring to the weak topology are preceded by $w$-.

\item We say that \emph{$\{A_{n} \colon n \in \nat\}$ converges to
    $A$ weakly polynomially}, $A = pw$-$\lim_{n \in \nat} A_{n}$
    in notation, if for every $k \in \nat$, $w$-$\lim_{n \in \nat}
    A_{n}^{k} = A^{k}$. The topology corresponding to this notion
    of convergence is called the \emph{weak polynomial topology}.
    Topological notions referring to the weak polynomial topology
    are preceded by $pw$-.

\item We say that \emph{$\{A_{n} \colon n \in \nat\}$ converges to $A$
    strongly}, $A = s$-$\lim_{n \in \nat} A_{n}$ in notation, if for
  every $x\in H$, $\lim_{n \in \nat} A_{n}x = Ax$.
The topology corresponding to this notion of convergence is called
the \emph{strong topology}. Topological notions referring to the
strong topology are preceded by $s$-.

\item We say that \emph{$\{A_{n} \colon n \in \nat\}$ converges to
    $A$ in the strong-star sense}, $A = s^{\star}$-$\lim_{n \in
    \nat} A_{n}$ in notation, if $s$-$\lim_{n \in \nat} A_{n} = A$
    and $s$-$\lim_{n \in \nat} A^{\star}_{n} = A^{\star}$. The
    topology corresponding to this notion of convergence is called
    the \emph{strong-star topology}. Topological notions referring
    to the strong-star topology are preceded by $s^{\star}$-. \een
\end{defin}

\begin{prop}\label{topolporp}(\cite[Section 2 p.\ 67]{T}, \cite[Section 7.f p.\ 121]{W}) Let $\{e_{i} \colon i \in \nat\} \ss H$ be an orthonormal basis.

\ben
\item For every $A,B \in B(H)$, set $$d_{w}(A,B) = \sum_{i,j \in
    \nat} 2^{-i-j} |\sp{Ae_{i},e_{j}} - \sp{Be_{i},e_{j}}|.$$ Then
    $d_{w}$ is a complete separable metric on $C(H)$ which
    generates the weak topology.

\item For every $A,B \in B(H)$, set $$d_{s}(A,B) = \sum_{i \in
    \nat} 2^{-i} \|Ae_{i} - Be_{i}\|.$$ Then $d_{s}$ is a complete
    separable metric on $C(H)$ which generates the strong
    topology.

\een
\end{prop}
The strong-star topology on uniformly bounded sets is generated by
the metric $d_{s^\star}(A,B)=d_s(A,B)+d_s(A^{\star},B^{\star})$.
It is easy to see that $d_{s^\star}$ is a complete metric on
$C(H)$. We will see in Section \ref{sspw} that $C(H)$ endowed with
the weak polynomial topology is a Polish space as well.
%We mention without proof that $C(H)$ endowed with the weak polynomial topology
%is \emph{not} a Polish space; nevertheless, as we will see in Section \ref{sspw},
%it is a Baire space.

Note that the weak, weak polynomial, strong, strong-star, and norm
topologies all refine the preceding topologies in this list. We
also need the following.

\begin{prop}\label{unitua}(\cite[Remark 4.10 p.\ 84]{T}) On $U(H)$ the weak, weak polynomial, strong, and strong-star topologies coincide. With this topology, $U(H)$ is a Polish space.
\end{prop}

%\subsection{The weak topology}\label{ssweak}
\section{The weak topology}\label{ssweak}

In the weak topology, the theories of typical properties of
contractions and of unitary operators coincide.

\begin{thm}\label{typunit}(\cite[Theorem 2.2 p.\ 2]{E0}) A $w$-typical contraction is unitary.
\end{thm}

Equivalently, Theorem \ref{typunit} says that $U(H)$ is a
$w$-co-meager subset of $C(H)$.
By Proposition \ref{unitua},
$U(H)$ endowed with the weak topology is also a Polish
space.
Hence the notions related to Baire category make sense relative to $U(H)$, and a set $A \ss C(H)$ is $w$-co-meager in $C(H)$ if and only if $A \cap U(H)$ is $w$-co-meager in $U(H)$.

Unitary operators are well-understood. E.g.\ the theory of spectral measures allows a detailed description of the spectral properties and conjugacy classes of unitary operators. We refer to \cite[Chapter XI.4 p.\ 306]{Y} and \cite[Chapter 2 p.\ 17]{N} for an introduction to spectral measures. The following proposition briefly summarizes how the spectral properties of $w$-typical contractions can be obtained from the theory of spectral measures. We set $S^{1} = \{\lambda \in \cmpl \colon |\lambda| =1\}$. For a measure $\mu$ on $S^{1}$, $\supp \mu$ denotes the closed support of $\mu$. For measures $\mu, \nu$ on $S^{1}$, we write $\mu \bot \nu$ if $\mu$ and $\nu$ are mutually singular.

\begin{prop}\label{weakspectrum} A $w$-typical contraction $U$ satisfies $C_{\sigma}(U) = S^{1}$.
\end{prop}
\prf Let $U$ be a $w$-typical contraction. By Theorem \ref{typunit}, $U$ is unitary. Let $\mu_{U}$ denote the maximal spectral type of $U$ (see \cite[Section 8.22 p.\ 55]{N}), which is a Borel measure on $S^{1}$. E.g.\ by the uniqueness of the measure class of $\mu_{U}$, $\lambda \in \sigma(U)$ if and only if $\lambda \in \supp \mu_{U}$, and $\lambda \in P_{\sigma}(U)$ if and only if $\lambda$ is an atom of $\mu_{U}$.

By \cite[8.25 Theorem (a) p.\ 56]{N}, $$\{U \in U(H) \colon P_{\sigma}(U) = \es\} =  \{U \in U(H) \colon \mu_{U} \textrm{ is atomless}\}$$ is a $w$-co-meager set in $U(H)$. By a similar argument, using \cite[7.7 Corollary p.\ 46]{N}, we obtain that $$\{U \in U(H) \colon \sigma(U) = S^{1}\} =  \{U \in U(H) \colon \supp \mu_{U} = S^{1}\}$$ is also a $w$-co-meager set in $U(H)$. Since $\sigma(U) = P_{\sigma}(U) \cup C_{\sigma}(U) \cup R_{\sigma}(U)$ and $R_{\sigma}(U) = \es$ for every unitary operator, we conclude that a $w$-typical contraction $U$ satisfies $C_{\sigma}(U) = S^{1}$, as required.$\bs$

\medskip

By analogous applications of spectral measures, one can isolate numerous additional $w$-typical properties of contractions (see e.g.\ \cite[8.25 Theorem p.\ 56]{N}). We refer to e.g.~\cite[Section IV.3]{E} for asymptotic properties of $w$-typical
contractions, and mention that P.~Zorin \cite{Pav} showed recently
that
a $w$-typical contraction admits a fixed cyclic vector.
Here we restrict ourselves to pointing out that despite the abundance of $w$-typical properties, typical contractions are not unitarily equivalent.

\begin{prop}\label{weakeq} For every $U \in C(H)$, $\mc{O}(U)$ is $w$-meager in $C(H)$. In particular, for a $w$-typical pair of contractions $(U_{1}, U_{2} )\in C(H) \times C(H)$ we have that $U_{1}$ and $U_{2}$ are not unitarily equivalent.
\end{prop}
\prf Let $U \in C(H)$ be arbitrary. By Theorem \ref{typunit}, the statement follows if $U \notin U(H)$. So for the rest of the argument, we can assume $U \in U(H)$. By the conjugacy invariance of the measure class of the maximal spectral type, for every $V \in U(H)$ we have that $\mu_{U}$ and $\mu_{VUV^{-1}}$ are mutually absolutely continuous. By \cite[8.25 Theorem (b) p.\ 56]{N}, for every measure $\nu$ on $S^{1}$, the set $\{V \in U(H) \colon \mu_{V} \bot \nu\}$ is $w$-co-meager  in $U(H)$.  Thus $\mc{O}(U)$ is $w$-meager in $U(H)$ and so in $C(H)$, as well.

Finally consider the set
$$
E= \{(U_{1}, U_{2} )\in C(H) \times C(H) \colon U_{1}, U_{2} \textrm{
  are unitarily equivalent}\}.
$$
The set $E$ is clearly analytic hence has the Baire property (see e.g.\ \cite[(29.14) Corollary p.\ 229]{K}). As we have seen above, for every $U_{1} \in C(H)$ we have $\{U_{2 }\in C(H) \colon (U_{1}, U_{2} )\in E\}$ is $w$-meager in $C(H)$. So by the Kuratowski-Ulam theorem (see e.g.\ \cite[(29.14) Corollary p.\ 229]{K}) we get that $E$ is $w$-meager in $C(H) \times C(H)$. This completes the proof.$\bs$

\medskip

%\subsection{The weak polynomial topology}\label{sspw}
\section{The weak polynomial topology}\label{sspw}

The main result of this section is the following.

\begin{thm}\label{pwBaire} The set $C(H)$ endowed with the weak
polynomial topology is a Polish space. Moreover, $U(H)$ is a
$pw$-co-meager $pw$-$G_{\delta}$ subset of $C(H)$.
\end{thm}

This result immediately implies that the theories of $pw$-typical
and $w$-typical properties of contractions coincide.
\begin{cor}\label{pwBaire_cor} A set $\mc{C} \ss C(H)$ is $pw$-co-meager in $C(H)$
if and only if $\mc{C} \cap U(H)$ is $w$-co-meager in $U(H)$.
In particular, a property $\Phi$ of contractions is $pw$-typical if and only if
$\Phi$ is $w$-typical.
\end{cor}
\prf By Theorem \ref{pwBaire}, $U(H)$ is a $pw$-co-meager subset
of the Polish space $C(H)$. By Proposition \ref{unitua}, the weak
and weak polynomial topologies coincide on $U(H)$ and in these
topologies $U(H)$ is a Polish space. So the notions related to
Baire category make sense relative to $U(H)$. By $U(H) \ss C(H)$
being $pw$-co-meager, a set $\mc{M} \ss U(H)$ is $pw$-meager in
$U(H)$ if and only if $\mc{M} \ss C(H)$ is $pw$-meager in $C(H)$.

We obtained that $\mc{C} \ss C(H)$ is $pw$-co-meager in $C(H)$ if
and only if $\mc{C} \cap U(H) \ss C(H)$ is $pw$-co-meager in
$C(H)$. This is equivalent to $\mc{C} \cap U(H) \ss U(H)$ being
$pw$-co-meager in $U(H)$, which is the same as $\mc{C} \cap U(H)
\ss U(H)$ being $w$-co-meager in $U(H)$. Finally by Theorem
\ref{typunit}, this is equivalent to $\mc{C} \ss C(H)$ being
$w$-co-meager in $C(H)$. This completes the proof.$\bs$

\medskip

By Corollary \ref{pwBaire_cor}, for the $pw$-typical properties of
contractions one can refer to Section \ref{ssweak}.

\medskip

To prove the first part of Theorem \ref{pwBaire} we need the
following lemmas on the weak topology.

\begin{lem}\label{str-wk} Let $A, A_{n} \in B(H)$ $(n \in \nat)$ satisfy $A = w$-$\lim_{n \in \nat}
A_{n}$.
If for every $x \in S_{H}$ we have $\|Ax\| \geq \limsup_{n \in
\nat} \|A_{n}x\|$ then $A = s$-$\lim_{n \in \nat} A_{n}$.
\end{lem}
\prf For $x \in S_{H}$ we have
$$
 \|Ax - A_{n}x \|^{2} = \|Ax\|^{2} + \|A_{n}x\|^2 - 2\re \sp{Ax,A_{n}x},
$$
so
\begin{multline}
   \notag \TS 0 \leq \limsup_{n \in \nat} \|Ax - A_{n}x \|^{2} = \|Ax\|^{2} + \limsup_{n \in \nat} \|A_{n}x\|^2 - \\\TS
   2\lim_{n \in \nat} \re \sp{Ax,A_{n}x} \leq 2\|Ax\|^{2} - 2 \re \sp{Ax,Ax} = 0.
\end{multline}
Since $x\in S_{H}$ was arbitrary, $A = s$-$\lim_{n \in \nat}
A_{n}$ follows. $\bs$

\begin{lem}\label{vawopn} Let $n > 0$, $\{x_{i} \colon i < n\} \ss S_{H}$ and $B \in C(H)$ be arbitrary.
Then for every $ \ve > 0$ there exists a $w$-open set $W \ss C(H)$
such that $B \in W$ and for every $A \in W$ we have $\|Ax_{i}\|
\geq \|Bx_{i}\| - \ve$ $(i < n)$.
\end{lem}

\prf We prove the statement for $n=1$ only; for $n > 1$ the
required set $W$ can be obtained by intersecting the sets $W_{i}$
satisfying the conditions of the lemma for each $x_{i}$ $(i < n)$
separately.

If $Bx_{0}=0$ the statement is trivial; so we can assume
$Bx_{0}\neq 0$. Consider the set
$$
 W = \{A \in C(H) \colon |\sp{Ax_{0},Bx_0}|> \|Bx_0\|^2-\|Bx_0\|\ve \};
$$
then $W$ is $w$-open and $B \in W$. For every $A \in W$ we have
$$
  \|Ax_{0}\|\geq \frac{|\sp{Ax_{0},Bx_0}|}{\|Bx_0\|} \geq \|Bx_0\| - \ve
$$
as required. $\bs$

\textbf{Proof of Theorem \ref{pwBaire}.} By \cite[Theorem
2.2]{E0}, $U(H)$ is a $w$-co-meager $w$-$G_{\delta}$ subset of
$C(H)$. Since the weak polynomial topology is finer than the weak
topology, $U(H) \ss C(H)$ is $pw$-$G_{\delta}$, as well. By
\cite[Theorem 1 p.\ 1142]{P}, $U(H)$ is $pw$-dense in $C(H)$, thus
$U(H)$ is $pw$-co-meager in $C(H)$.

To prove that $C(H)$ is Polish in the weak polynomial topology we
use \cite[(8.18) Theorem p.\ 45]{K} stating that a non-empty,
second countable topological space is Polish if and only if it is
$T_{1}$, regular and strong Choquet (for the definition, see \cite[(8.14) Definition p.\ 44.]{K}). Since $(C(H),pw)$ is a metric
space, it is second countable, $T_{1}$ and regular. So in order to
show it is Polish, we have to prove it is strong Choquet, i.e.\
that player $II$ has a winning strategy in the strong Choquet
game.

We define the strategy for player $II$ as follows. As in  Proposition \ref{topolporp}, let $d_w$
denote the \emph{complete} metric on $C(H)$ which induces the
\emph{weak} topology on $C(H)$.
Let $\{x_{n} \colon n \in \nat\}$ be a dense subset of $S_{H}$.
Suppose the $n^{\mathrm{th}}$ move of player $I$ is
$(A_{n},U_{n})$, where $A_{n} \in U_{n} \ss C(H)$ and $U_{n}$ is a
$pw$-open set. Let $W_{n}$ be a $w$-open set with the following
properties: \ben
\item\label{balls} $A_{n} \in W_{n} \ss \cl_{w}(W_{n}) \ss \{A
    \colon
    d_w(A,A_{n}) < 1/(n+1)\}$;
\item\label{decrease} $\cl_{w}(W_{n})\ss W_{n-1}$;
\item\label{norm-ineq} for every $A \in W_{n}$ we have $\|Ax_{i}\|
    \geq
    \|A_{n}x_{i}\| - 1/(n+1)~(i \leq n);$
\een such a set exists by Lemma \ref{vawopn}.
    Then let player II respond by playing a
    $pw$-open set $V_{n} \ss U_{n}$ such that \ben
    \setcounter{enumi}{3}
\item\label{inclusion} $A_{n} \in V_{n} \ss \cl_{pw}(V_{n}) \ss U_{n} \cap
    W_{n}$.
    \een

We show that this strategy is winning for player $II$. Let
$\{(A_{n},U_{n}),V_{n} \colon n \in \nat\}$ be a run in the game
in which player $II$ follows the above strategy. By $A_{n+1} \in
V_{n} \ss W_{n}$ $(n \in \nat)$ and
conditions \ref{balls} and \ref{decrease}, $A_{n}$ is weakly
convergent, say $A = w$-$\lim_{n \in \nat} A_{n}$. Again by
conditions \ref{balls} and \ref{decrease}, $A \in W_{n}$ $(n \in
\nat)$; thus by condition \ref{norm-ineq}, $\|Ax_{i}\| \geq
\limsup_{n \in \nat} \|A_{n}x_{i}\|$ $(i \in \nat)$. Since
$\{x_{n} \colon n \in \nat\} \ss S_{H}$ is dense and $A,A_{n}$ $(n
\in \nat)$ are contractions, we get $\|Ax\| \geq \limsup_{n \in
\nat} \|A_{n}x\|$ $(x \in S_{H})$. Thus by Lemma \ref{str-wk}, we
have $A = s$-$\lim_{n \in \nat} A_{n}$; in particular, $A =
pw$-$\lim_{n \in \nat} A_{n}$. Since $A \in \cl_{pw}(V_{n})\ss
V_{n-1}$ for each $n$ by condition \ref{inclusion}, we get $A \in
\bigcap_{n \in \nat} V_{n}$. This shows that the strategy is
winning for player $II$, and finishes the proof. $\bs$

\smallskip

Note that, by the same argument as in the above proof, $C(H)$ is
Polish for any metrizable topology on $C(H)$ which is finer than
the weak topology and coarser than the strong one. In addition,
for such topology $U(H)$ is a co-meager $G_\delta$ subset of
$C(H)$ if and only if $U(H)$ is dense in $C(H)$ in this topology.

%\subsection{The strong topology}\label{ssstrong}
\section{The strong topology}\label{ssstrong}

Probably the most surprising observation in the present paper is that
some typical properties of contractions in the strong and weak
topologies are completely different. While typical contractions in the
weak topology are not unitarily equivalent, typical contractions in
the strong topology are unitarily equivalent to an infinite
dimensional backward unilateral shift operator, hence the
investigation of $s$-typical properties of contractions is reduced to the study of one particular operator. We
introduce this operator in the following definition.
Since infinite dimensional complex separable Hilbert spaces are isometrically isomorphic, we can restrict ourselves to the study of a particular one.

\begin{defin}\label{Shift}\rm Set $H=\ell^{2}(\nat \times \nat)$ and denote the canonical orthonormal basis of $H$ by $\{e_{i}(n)   \colon i,n \in \nat \}$. We define the \emph{infinite dimensional backward unilateral shift operator} $S \in C(H)$ by $Se_{0}(n) = 0$ and $Se_{i+1}(n) = e_{i}(n)$ $(i,n \in \nat)$.
\end{defin}

The main result of this section is the following. Recall $\mc{O}(A) = \{UAU^{-1} \colon U \in U(H)\}$.

\begin{thm}\label{ueshift} The set $\mc{O}(S)$ is an $s$-co-meager subset of $C(H)$.
\end{thm}

By Theorem \ref{ueshift}, a property of contractions is $s$-typical if and only if $S$ has this property. In the following corollary, we recall only the properties of $S$ we usually concern in this paper.

\begin{cor}\label{ueshift_cor} An $s$-typical contraction $A$ satisfies that
\ben
\item \label{ueshift_cor_1} $P_{\sigma}(A) = \{\lambda \in \cmpl \colon |\lambda | < 1\}$, and for every $\lambda \in P_{\sigma}(A)$ we have $\dim \ker (\lambda \cdot \Id - A) = \infty$;

\item \label{ueshift_cor_2} $C_{\sigma}(A) = S^{1}$;

\item \label{ueshift_cor_3} $A$ can be embedded into a strongly
  continuous semigroup.
\een
Moreover, typical contractions are unitarily equivalent.
\end{cor}
\prf Statements \ref{ueshift_cor_1} and \ref{ueshift_cor_2} follow from Theorem \ref{ueshift} and the results of \cite{S}, while
\ref{ueshift_cor_3} is a corollary of \cite[Proposition 4.3]{E0}.
$\bs$

\medskip

 Our strategy to prove Theorem \ref{ueshift} is the following. The
 main observation is that for an $s$-typical contraction $A$, its
 \emph{adjoint} $A^{\star}$ is an isometry. Then by the Wold
 decomposition theorem (see e.g.\ \cite[Theorem 1.1 p.\ 3]{SzF}), $A$
 is unitarily equivalent to a direct sum of unitary and
backward unilateral shift operators, and the number of shifts in the direct sum depends on the dimension of $\ker A$. Since an $s$-typical contraction $A$ is strongly stable, the unitary part is trivial. So we complete the proof by showing $\dim \ker A = \infty$.

%\subsubsection{Elementary observations}
\subsection{Elementary observations}

 We collect here some elementary results we need later in our analysis.

 \begin{lem}\label{sum_norm} Let $x,y \in S_{H}$ satisfy $x \neq -y$. Set $\alpha =(2+2\cdot \re\sp{x,y})^{-1/2}$. Then $\|\alpha(x+y)\| = 1$.
\end{lem}
\prf We have $\|\alpha(x+y)\|^{2} = \alpha^2 (\|x\|^{2} + \|y\|^{2} + 2\cdot \re\sp{x,y}) = \alpha^2 (2 + 2\cdot \re\sp{x,y}) = 1$, as required.$\bs$

\begin{lem}\label{egynorm} Let $A \in C(H)$ and let $(b_{n})_{n \in \nat} \ss S_{H}$, $z \in S_{H}$ satisfy $\lim _{n \in \nat} A b_{n} = z$. Then $(b_{n})_{n \in \nat}$ is convergent.
\end{lem}
\prf It is enough to prove that for every $\ve > 0$ there is an $N \in
\nat$ such that for every $n,m \geq N$ we have $\|b_{n} - b_{m}\|^2
\leq \ve$. So let $\ve > 0$ be arbitrary. For every $n,m$ sufficiently
large we have $b_n \neq -b_m$, so $\alpha_{n,m}=(2+2\cdot
\re\sp{b_n,b_m})^{-1/2}$ is defined. By Lemma \ref{sum_norm} we have
$\|\alpha_{n,m}(b_n + b_m)\| = 1$.
So since $A$ is a contraction,
\begin{multline} \label{m2}
 \alpha_{n,m} \cdot (2  -\|Ab_n-z\| -
  \|Ab_m-z\|)\leq \\ \alpha_{n,m}\cdot \|2 \cdot z +(Ab_n-z) + (Ab_m-z)\|
  = \\ \alpha_{n,m}\cdot \|Ab_n + Ab_m\| = \|\alpha_{n,m}\cdot A(b_n +
  b_m)\| \leq 1.
\end{multline}

Let $N \in \nat$ be such that for every $n \geq N$ we have $\|Ab_n - z\| \leq \ve/8$. Then by (\ref{m2}), for every $n,m \geq N$ we get $\alpha_{n,m} \leq 1/(2-2\ve/8)$, i.e.\ $2+2 \cdot \re \sp{b_n, b_m} \geq 4 \cdot (1-\ve/8)^2$. Thus $$\|b_n - b_m\|^2 = 2-2\cdot \re\sp{b_n,b_m} \leq 4-4\cdot (1-\ve/8)^2 \leq \ve,$$ as required.$\bs$

 \medskip

 The following lemma will be helpful to show that the kernel of a typical contraction is infinite dimensional.

\begin{lem}\label{ortho} Let $n \in \nat \sm \{0\}$ and let $\{e_{i}\colon i < n\} \ss H$ be an orthonormal family. Let $\{f_{i}\colon i < n\} \ss H$ satisfy $\|f_{i} - e_{i}\| < 1/n$ $(i < n)$. Then $\{f_{i}\colon i < n\}$ are linearly independent.
\end{lem}
\prf Let $\alpha_{i} \in \cmpl$ $(i < n)$ be arbitrary satisfying $\sum_{i < n}|\alpha_{i}| > 0$. We have \begin{multline} \notag \b\|\sum_{i < n} \alpha_{i}e_{i} - \sum_{i < n} \alpha_{i}f_{i}\j\|^{2} \leq \b(\sum_{i < n}|\alpha_{i}|\|e_{i} - f_{i}\|\j)^{2} < \\ \b(\frac{1}{n} \sum_{i < n}|\alpha_{i}|\j)^{2} \leq \frac{1}{n} \sum_{i < n}|\alpha_{i}|^{2} = \frac{1}{n}\b\|\sum_{i < n} \alpha_{i}e_{i}\j\|^{2}.\end{multline} So by $\b\|\sum_{i < n} \alpha_{i}e_{i}\j\| > 0$ we have $\b\|\sum_{i < n} \alpha_{i}f_{i}\j\| > 0$, as required.$\bs$

 \medskip

Next we point out a trivial sufficient condition for the direct sum of contractions to be a contraction.

 \begin{lem}\label{sum} Let $V_{0}, V_{1} \leq H$ be subspaces satisfying $V_{0} \bot V_{1}$. Let $A_{i} \colon V_{i} \rar H$ $(i < 2)$ be contractive linear operators such that $A_{0}[V_{0}] \bot A_{1}[V_{1}]$. Then $A \colon \span{V_{0},V_{1}} \rar H$, $$A(\alpha v_{0} + \beta v_{1}) = \alpha A_{0}v_{0} + \beta A_{1}v_{1}~(\alpha, \beta  \in \cmpl)$$ is also contractive.
\end{lem}
\prf Let $v \in \span{V_{0},V_{1}}$, $\|v\|=1$ be arbitrary. Then there are $\alpha, \beta  \in \cmpl$ with $|\alpha|^{2} + |\beta|^{2} = 1$ and $v_{i} \in S_{V_{i}}$ $(i < 2)$ such that $v = \alpha v_{0} + \beta v_{1}$. We have \begin{multline} \notag \|Av\|^{2} = \|A(\alpha v_{0} + \beta v_{1})\|^{2} = \|\alpha A_{0}v_{0} + \beta A_{1}v_{1}\|^{2} = \\ |\alpha|^{2}\|A_{0}v_{0}\|^{2} + |\beta|^{2}\|A_{1}v_{1}\|^{2} \leq |\alpha|^{2} + |\beta|^{2} = 1,\end{multline} so the proof is complete.$\bs$

\medskip

 Finally we point out that strongly stable contractions form a $s$-co-meager $s$-$G_{\delta}$ subset of $C(H)$.

 \begin{defin}\label{sstabledef}\rm A contraction $A$ is
   \emph{strongly stable} if $s$-$\lim_{n \in \nat}A^n=0$, i.e.,
 for every $x \in S_{H}$ and $\ve > 0$ there is an $n \in \nat$ such that $\|A^{n}x\| < \ve$. The set of strongly stable contractions is denoted by $\mc{S}$.
 \end{defin}

 \begin{lem}\label{stable} The set of strongly stable contractions is an $s$-co-meager $s$-$G_\delta$ subset of $C(H)$.
\end{lem}
\prf
By $A=\lim_{n \in \nat}(1-2^{-n})A$ $(A \in C(H))$, the set of
contractions  $A$ satisfying $\|A\|<1$ is a norm dense and hence
an $s$-dense subset of $C(H)$. Since every such operator is strongly
stable, it remains to show that $\mc{S}$ is $s$-$G_\delta$. To this end, let $\{x_i \colon i \in \nat\}$ be a  dense subset of
$S_{H}$.  Note that for each $x\in H$ and $A \in \mc{S}$, the sequence $\|A^nx\|$ $(n \in \nat)$ is monotonically decreasing, so $\lim_{n\in\nat}\|A^nx\|=0$ is equivalent to $\inf_{n\in\nat}\|A^nx\|=0$. Thus
\Keq\label{Juj}
\mc{S} = \TS \bigcap_{i,j \in \nat} \bigcup_{n \in \nat} \{A \in C(H) \colon  \|A^n x_i\|< 2^{-j}\},
\Zeq which completes the proof.$\bs$

\medskip

We remark that  Lemma \ref{stable} holds in every separable Banach space. It is interesting to note the difference to the weak operator
topology in which the set of weakly stable contractions is $w$-meager
in $C(H)$ (see \cite[Theorem 4.3]{ES1}).

%\subsubsection{Mapping properties of $s$-typical contractions}
\subsection{Mapping properties of $s$-typical contractions}

The purpose of this section is to prove the following.

\begin{thm}\label{styp} Let $\mc{G}$ denote the set of contractive operators $A$
satisfying the following properties:
\begin{enumerate}
\item\label{styp_1} for every $y \in S_{H}$ there exists an $x \in
  S_{H}$ such that $Ax=y$, i.e., $S_{H} \ss A[S_{H}]$;
\item\label{styp_2} $\dim \Ker A=\infty$.
\end{enumerate}
Then  $\mc{G}$ is an $s$-co-meager subset of $C(H)$.
\end{thm}

To prove Theorem \ref{styp}, we need a geometric lemma saying that every contraction defined on a finite dimensional subspace of $H$ can be extended to a contraction which is surjective in a very strong sense.

\begin{lem}\label{ext} Let $V \leq H$ be a finite dimensional subspace and let
$A \colon V \rar H$ be a contractive linear operator. Let $W = A[V]$ and let $Y \leq H$ be an arbitrary subspace satisfying $Y \bot W$. Then for every $X \leq H$ satisfying $X \bot V$ and $\dim X = \dim W + \dim Y$ there exists a contractive linear operator $\tilde A \colon \span{V,X} \rar H$ such that $\tilde A|_{V} = A$ and $\tilde A[B_{\span{V,X}}] = B_{\span{W,Y}}$.
\end{lem}
\prf We handle first the special $Y = \{0\}$ case. Let $\dim V=n$ and $ \dim W = m$.
The Gram-Schmidt orthogonalization theorem states that there exists
an orthonormal base $\{v_{i} \colon  i < n\} \ss V$ such that
$\{Av_{i} \colon i < n\} \ss W$ are pairwise orthogonal, and
\ben[(i)]
\item\label{fe1} $\|Av_{0}\| = \max\{\|Av\| \colon v \in S_{V}\}$;
\item\label{fe2} for every $j < n-1$, $\|Av_{j+1}\| = \max\{ \|Av\| \colon v \in S_{V} \cap \span{ v_{i} \colon i \leq j} ^{\bot}\}$.
\een
By (\ref{fe1}) and (\ref{fe2}) we have
$\|Av_{i}\| > 0$ if and only if $i < m$. Set $w_{i} = Av_{i}/\|Av_{i}\|$ $(i < m)$.

 Let  $X \leq H$ satisfy $X \bot V$ and $\dim X = m$. Fix an orthonormal base $\{x_{i} \colon i < m\}$ in $X$ and define $\tilde A \colon \span{V,X} \rar H$ such that $\tilde A|_{V} = A$ and $\tilde A x_{i} =  \sqrt{1-\|Av_{i}\|^{2}}w_{i}$ $(i < m)$.

First we show $\tilde A$ is a contraction. Let $u \in \span{V,X}$ satisfy $\|u\| = 1$. Then $u = \sum_{i < n} \alpha_{i} v_{i} + \sum _{i < m} \beta_{i} x_{i}$ where $\alpha_{i} \in \cmpl$ $(i < n)$,  $\beta_{i} \in \cmpl$ $(i < m)$ satisfy $\sum_{i < n} |\alpha_{i}|^{2} + \sum_{i < m} |\beta_{i}|^{2} = 1$. Then \begin{multline}\notag \|\tilde A u\|^{2} = \b\|\tilde A \b(\sum_{i < n} \alpha_{i} v_{i} + \sum _{i < m} \beta_{i} x_{i} \j) \j\|^{2} = \\   \b\|\sum _{i < m} (\alpha_{i}\|Av_{i}\|w_{i} + \beta_{i}\sqrt{1-\|Av_{i}\|^{2}}w_{i}) \j\|^{2} \leq \sum _{i < m} (|\alpha_{i}|\|Av_{i}\| + |\beta_{i}|\sqrt{1-\|Av_{i}\|^{2}})^{2}\end{multline}
By the Cauchy-Schwarz inequality, for every $0 \leq p,q,r \leq 1$ we have $(pr + q\sqrt{1-r^2})^{2} \leq p^{2} + q^{2}$. So $\|\tilde A u\| \leq \sum _{i < m} |\alpha_{i}|^2 + |\beta_{i}|^{2} \leq 1$, as required.

Next we show that for every $y \in  B_{W}$ there is an $x \in B_{\span{V,X}}$ such that $\tilde A x = y$. Let $y = \sum _{i < m} \beta_{i} w_{i}$ where $\beta_{i} \in \cmpl$ $(i < m)$ satisfy $\sum_{i < m} |\beta_{i}|^{2} \leq 1$. Let $$x = \sum_{i < m} \beta_{i}(\|Av_{i}\|v_{i} + \sqrt{1-\|Av_{i}\|^{2}}x_{i}).$$ Then $$\|x\| = \sum _{i <m} \b(\beta_{i}^{2}\|Av_{i}\|^2 + \beta_{i}^{2}\sqrt{1-\|Av_{i}\|^{2}}^2\j) \leq 1,$$ and \begin{multline}\notag\tilde A x = \sum_{i < m} \beta_{i}(\|Av_{i}\|Av_{i} + \sqrt{1-\|Av_{i}\|^{2}} \tilde A x_{i}) = \\  \sum_{i < m} \beta_{i}(\|Av_{i}\|\|Av_{i}\|w_{i} + \sqrt{1-\|Av_{i}\|^{2}} \sqrt{1-\|Av_{i}\|^{2}} w_{i}) = y,\end{multline} as required. This completes the proof of the special $Y = \{0\}$ case.

In the general case write $X = X_{0} \oplus X_{1}$ where $X_{0} \bot X_{1}$ and $\dim X_{0} = \dim W$, $\dim X_{1} = \dim Y$. By the special case above, there is a contraction $\tilde A \colon \span{V,X_{0}} \rar H$ such that $\tilde A|_{V} = A$ and $\tilde A[B_{\span{V,X_{0}}}] = B_{W}$. Extend further $\tilde A$ by setting $\tilde A|_{X_{1}} \colon X_{1} \rar Y$ be any isometric isomorphism. By Lemma \ref{sum}, $\tilde A$ is a contraction which clearly satisfies  $\tilde A[B_{\span{V,X}}] = B_{\span{W,Y}}$. This completes the proof.$\bs$

\medskip

From Lemma \ref{ext}, we immediately get the following.

\begin{prop}\label{dense} The set of
contractive operators $A$ such that
for every
$y \in S_{H}$
there exists an $x \in  S_{H}$ such that
$Ax=y$ is an $s$-co-meager subset of $C(H)$.
\end{prop}
\prf
Let $$\mc{M} = \{ A \in C(H) \colon \forall \ve > 0, y \in S_{H} ~\exists x \in  S_{H} ~(\|y - Ax\| < \ve)\}.$$ First we show that $\mc{M}$ is an $s$-dense $s$-$G_\delta$ subset of $C(H)$.

Fix $y \in S_{H}$ and $\ve > 0$. The set
$$
 C(y,\ve) = \{A \in C(H) \colon  \exists x \in S_{H} ~(\|y - Ax\| < \ve)\}
$$
is $s$-open. We show that it is $s$-dense.

Let $U \ss C(H)$ be any non-empty $s$-open set. By passing to a
subset, we can assume $U=\{A \in C(H) \colon \|y_{i} - Ax_{i}\| <
\ve_{i} ~(i \in I)\}$ where $x_{i},y_{i} \in H$, $\ve_{i} > 0$ $(i
\in I)$ and $I$ is finite.  Let $V = \span{ x_{i}
\colon i \in I}$ and take an arbitrary $A \in U$. By restricting $A$
to $V$ we can assume $A|_{V^{\bot}} = 0$. Set $W = A[V]$.

Let $Y \leq H$ be an at most one dimensional subspace such that $Y \bot W$ and $y \in \span{W,Y}$. Let $X \leq H$ be a $\dim W + \dim Y$ dimensional subspace such that $X \bot V$. By Lemma \ref{ext}, there exists a contraction $\tilde A \colon \span{V,X} \rar H$ such that $\tilde A|_{V} = A|_{V}$ and $\tilde A[B_{\span{V,X}}] = B_{\span{W,Y}}$. In particular, there is an $x \in  B_{\span{V,X}}$ such that $\tilde A x = y$. Since $\tilde A$ is a contraction and $\|y\|=1$, we get $x \in S_{H}$. Extend further $\tilde A$ by setting $\tilde A|_{\span{V,X}^{\bot}} = 0$. Then $\tilde A \in U \cap C(y,\ve)$; i.e.\ we concluded that $C(y,\ve)$ is $s$-dense.

Let $D \ss S_{H}$ be a countable dense set. We have $$
\mc{M} =   \TS \bigcap\{C(y,2^{-n}) \colon y \in D,~n \in \nat\},
$$
so by the Baire category theorem, $\mc{M}$ is an $s$-dense $s$-$G_{\delta}$
subset of $C(H)$.

It now suffices to show that every $A\in \mc{M}$ satisfies $S_{H} \ss A[S_{H}]$. To this end, let $A\in \mc{M}$ and $z \in S_H$ be arbitrary. By the definition
of $\mc{M}$, there is a sequence $(b_{n})_{n \in \nat} \ss S_{H}$ such that
$\lim _{n \in \nat} A b_{n} = z$.
By Proposition \ref{egynorm}, $(b_{n})_{n \in \nat}$ is convergent,
say $\lim _{n \in \nat}  b_{n} = x$.
Then $Ax=z$, which completes the proof.
$\bs$

\medskip

To prove that a typical contraction has infinite dimensional kernel, we need a lemma showing that a typical contraction approximates the zero operator on arbitrarily large finite dimensional subspaces.

\begin{lem}\label{finite} The set of contractive operators $A$ such
  that for every $n \in \nat$ and $\ve > 0$ there exists $Z \leq H$
  with $\dim Z \geq n$ and $\|A|_{Z}\| < \ve$ is
an $s$-dense $s$-$G_\delta$
subset of $C(H)$.
\end{lem}
\prf For every $n \in \nat$ and $\ve > 0$, the set $$ C(n,\ve) = \{A \in C(H) \colon \exists Z\leq  H~(\dim Z \geq n,~\|A|_{Z}\| < \ve)\}$$ is $s$-open. We show that $C(n,\ve)$ $(n \in \nat,\ve > 0)$ are $s$-dense.

Fix arbitrary $n \in \nat$ and $\ve > 0$. Let $U \ss C(H)$ be any non-empty $s$-open set. By passing to a
subset, we can assume $U=\{A \in C(H) \colon \|y_{i} - Ax_{i}\| <
\ve_{i} ~(i \in I)\}$ where $x_{i},y_{i} \in H$, $\ve_{i} > 0$ $(i
\in I)$ and $I$ is finite. Let  $A \in U$ be arbitrary. Let $V = \span{ x_{i}
\colon i \in I}$, and define $B \in C(H)$ by $B|_{V} = A|_{V}$,
$B|_{V^{\bot}} = 0$. Since $\dim V^{\bot} = \infty$, we obtained $B \in C(n,\ve) \cap U$; i.e.\ we concluded $C(n,\ve)$ is $s$-dense.

The set of contractive operators $A$ such that for every $n \in \nat$ and $\ve > 0$ there exists $Z \leq H$ with $\dim Z \geq n$ and $\|A|_{Z}\| < \ve$ is $$\TS \bigcap\{C(n,2^{-m}) \colon n, m \in \nat\}.$$ So by the Baire category theorem, this is an $s$-dense $s$-$G_{\delta}$ subset of $C(H)$, which completes the proof.$\bs$

\medskip

We are ready to prove the second part of Theorem \ref{styp}.

\begin{prop}\label{kern} The set of
contractions $A$ which satisfy $\dim \ker A = \infty$ is an $s$-co-meager subset of $C(H)$.
\end{prop}
\prf By Lemma \ref{finite} and Proposition \ref{dense}, the set of
contractions $A$ which satisfy that
\ben
\item\label{e1} for every $n \in \nat$ and $\ve > 0$ there exists $Z \leq H$ with $\dim Z \geq n$ and $\|A|_{Z}\| < \ve/n$;
\item\label{e2} for every $y \in S_{H}$
there exists an $x \in  S_{H}$ such that
$Ax=y$;
\een
 is an $s$-co-meager subset of $C(H)$. We show that every member $A$ of this set satisfies $\dim \ker A = \infty$; this will complete the proof.

Fix an arbitrary $n \in \nat \sm \{0\}$. Set $\ve = 1$ and let $Z \leq H$ satisfy \ref{e1} for this $n$ and $\ve$. Let $\{e_{i}\colon i < n\} \ss Z$ be an orthonormal family. We find $f_{i} \in H$ $(i < n)$ such that $\|f_{i} - e_{i}\| < 1/n$ $(i < n)$ and $Af_{i} = 0$. Then by Lemma \ref{ortho}, $\{f_{i}\colon i < n\}$ are linearly independent. Then $\dim \ker A \geq n$, and since $n$ was arbitrary, we concluded $\dim \ker A= \infty$.

Fix $i<n$; we define $f_i$ as follows. If $Ae_{i} = 0$, set $f_{i} = e_{i}$. Else observe that by \ref{e2}, there
exists $x_i\in S_H$ with $Ax_i= Ae_i / \|Ae_i\|$. Define now
$$
  f_i=e_i-\|Ae_i\|x_i
.$$
Then $\|f_i-e_i\|=\|Ae_i\|<1/n$ and we have
$
  Af_i=Ae_i-\|Ae_i\| Ax_i=0$,
so the construction is complete.$\bs$

\medskip

\textbf{Proof of Theorem \ref{styp}.} The statement follows from Proposition \ref{dense} and Proposition \ref{kern}.$\bs$

%\subsubsection{Unitary equivalence of $s$-typical contractions to a shift}
\subsection{Unitary equivalence of $s$-typical contractions to a shift}

Operators in the set $\mc{G}$ introduced in Theorem \ref{styp} have the following properties.

\begin{prop}\label{AA} Let $A \in \mc{G}$. Then $AA^{\star} = \Id$ and
  $A^{\star}A$ is the projection onto $\Ran A^{\star}$, which is an
  infinite dimensional and infinite co-dimensional subspace of $H$. In
  particular, $A^{\star}$ is an isometry, hence $A$ is a co-isometry
  and in addition, $A$ is an isometry on $(\Ker A)^{\bot}$.
\end{prop}
\prf Let $\{e_{ i} \colon i \in \nat\}$ be an orthonormal basis of $H$.  By the definition of $\mc{G}$, for every $i \in \nat$ there is an $a_{i} \in S_{H}$ such that $Aa_{i} = e_{i}$. Note that for every $i \in \nat$, $1 = \sp{e_{i},e_{i}} = \sp{Aa_{i}, Aa_{i}} = \sp{A^{\star} A a_{i}, a_{i}}$.  By $A$, $A^{\star}$ being contractions, this is possible only if $A^{\star} e_{i} = A^{\star}Aa_{i} = a_{i}$ $(i \in \nat)$, thus $AA^{\star}e_{i} = Aa_{i}=e_{i}$ $(i \in \nat)$. This proves that $AA^{\star} = \Id$, $\Ran A^{\star} = \span{A^{\star}e_{i} \colon i \in \nat} = \span{a_{i} \colon i \in \nat}$, and that $A^{\star}A$ is the projection onto $\Ran A^{\star}$.
Again by $A, A^{\star} \in C(H)$, this implies that $A^{\star}$ is an isometry, and $A$ is an isometry on $\Ran A^{\star}$.
By the definition of $\mc{G}$,
$\Ran A^{\star} = (\Ker A)^{\bot}$ is infinite dimensional and infinite co-dimensional. This completes the proof.
$\bs$

\medskip

\textbf{Proof of Theorem \ref{ueshift}.} By Lemma \ref{stable} and Theorem \ref{styp}, it is enough to show that every $A \in \mc{S} \cap \mc{G}$ is unitarily equivalent to the operator $S$ of Definition \ref{Shift}. By Proposition \ref{AA}, $A$ is a co-isometry. So by the Wold decomposition theorem (see e.g.\ \cite[Theorem 1.1 p.\ 3]{SzF}),  we have $H=H_u \oplus H_s$ such that
$A|_{H_u}$ is unitary and $A|_{H_s}$ is unitarily equivalent to the
backward unilateral shift operator on $l^2(\nat, \ker A)$, i.e.\ the Hilbert space of square summable $\nat \rar \ker A$ functions. By $A \in \mc{S}$ we have $H_u=\{0\}$. By Theorem \ref{styp}.\ref{styp_2}, $\dim \ker A = \infty$ hence $l^2(\nat, \ker A)$ is isometrically isomorphic to $\ell^{2}(\nat \times \nat)$. This completes the proof.$\bs$

%\subsubsection{The continuous case}\label{sstrong-cont}
\subsection{The continuous case}\label{sstrong-cont}

In this section we show that a typical strongly continuous contraction
semigroup on $H$ is unitarily
equivalent to an infinite dimensional backward unilateral shift semigroup.

Let $C^c(H)$ denote the set of contractive $C_0$-semigroups on $H$. Here we endow $C^c(H)$ with the topology induced by the uniform strong convergence on compact time intervals, i.e., by the
metric
$$
\TS
d_s^c(T(\cdot), S(\cdot)) = \sum_{j,n \in \nat} 2^{-(j+n)} \sup_{t\in[0,n]}\|T(t)e_j -S(t)e_j\|,
$$
where $\{e_j \colon j \in \nat\}$ is an orthonormal basis of $H$. With respect to this topology, $C^c(H)$ is a Polish space. With an abuse of notation, topological notions referring to this topology are also preceded by $s$-.

We again restrict ourselves without loss of generality to a particular
infinite dimensional complex separable Hilbert space and introduce on it the infinite dimensional backward unilateral shift semigroup.

\begin{defin}\rm
Set $H=L^2(\real_+^2)$. We define the \emph{infinite dimensional backward unilateral shift semigroup}
$S(\cdot)\in C^c(H)$ by
$$
  [S(t)f](s,w)=f(s+t,w) \quad (t,s,w\in \real_+,\ f\in H).
$$
We also set $\mc{O}(S(\cdot)) = \{(U^{-1}S(t)U)_{t\in \real_+} \colon U\in
U(H)\}\subseteq C^c(H)$.
\end{defin}

The following shows that an $s$-typical contraction semigroup is unitarily equivalent to $S(\cdot)$.

\begin{thm}\label{UE-cont}
The set $\mc{O}(S(\cdot))$ is an $s$-co-meager subset of $C^c(H)$.
\end{thm}

Note that for every $t > 0$, $S(t)$ is unitarily equivalent to the backward unilateral shift operator $S$ of Definition \ref{Shift}. So by Theorem \ref{UE-cont}, $s$-typical contractive $C_0$-semigroups $T(\cdot)$ satisfy that for every $t>0$, the operator $T(t)$ has the same properties as $S$. Moreover, $s$-typical contraction semigroups are unitarily equivalent.

To prove Theorem \ref{UE-cont}, we need to introduce the following not so well-known concept from
semigroup theory.
\begin{defin}\rm \label{CGN}
Let $T(\cdot)$ be a $C_0$-semigroup on $H$. Let the generator
$A$ of $T(\cdot)$ satisfy $1\in\rho(A)$. Then the operator $V \in B(H)$ defined by
$$
  V=(A+\Id)(A-\Id)^{-1}=\Id+2(A-\Id)^{-1}
$$
 is called the \emph{cogenerator} of $T(\cdot)$.
\end{defin}

The cogenerator is a bounded operator which determines the semigroup
uniquely. Moreover, it shares many properties of $T(\cdot)$ such as being contractive,
unitary, self-adjoint, normal, isometric, strongly stable etc.\ (see
\cite[Section III.8-9]{SzF} for details). We will use the following (see e.g.\ \cite[Theorem III.8.1]{SzF}).

\begin{lem} \label{eq:cogen} With the notation of Definition \ref{CGN},
$V$ is the cogenerator of a contraction semigroup if and only if $V$ is contractive and satisfies $1\notin P_\sigma(V)$.
\end{lem}

Our key observation is the following.

\begin{lem}\label{UCSO} Set $$\mc{V} = \{V \in C(H) \colon 1\notin P_\sigma(V)\},$$ endowed with the strong topology. Define $J \colon C^c(H)\rar \mathcal{V}$ by $J(T(\cdot)) = V$
where $V$ is the cogenerator of $T(\cdot)$ $(T(\cdot) \in C^{c}(H))$. Then $J$ is a homeomorphism.
\end{lem}
\prf The statement follows from Lemma \ref{eq:cogen} and the First Trotter--Kato
Theorem (see e.g.\ \cite[Theorem III.4.8]{EN}).$\bs$

\medskip

\textbf{Proof of Theorem \ref{UE-cont}.}
By Corollary \ref{ueshift_cor}, $\mc{V}$ is $s$-co-meager in $C(H)$. So with the notation of (\ref{Juj}) and Theorem \ref{styp}, by Lemma \ref{stable} and Theorem \ref{styp} we have that $\mc{S}\cap\mc{G}\cap \mc{V}$ is $s$-co-meager in $\mc{V}$. Hence the set $J^{-1}(\mc{S}\cap \mc{G}\cap \mc{V})$ is $s$-co-meager in $C^c(H)$. So it is enough to show that  $J^{-1}(\mc{S}\cap \mc{G}\cap
\mc{V}) \ss \mc{O}(S(\cdot))$.

Let $T(\cdot)\in J^{-1}(\mc{S}\cap \mc{G}\cap
\mc{V})$ be arbitrary. Let $V$ denote the cogenerator of $T(\cdot)$; then $V \in \mc{S}\cap \mc{G}\cap
\mc{V}$.
By Proposition \ref{AA}, $V$ is a strongly stable co-isometry. Since the cogenerator of
$T^\star(\cdot)$ is $V^\star$, and since the semigroup and the cogenerator
share strong stability and the isometric property (see
\cite[Theorem III.9.1]{SzF}), $T(\cdot)$ is a
strongly stable co-isometric contraction semigroup. By Wold's
decomposition for semigroups (see \cite[Theorem III.9.3]{SzF}), $T(\cdot)$ is
unitarily equivalent to the backward unilateral shift semigroup on $L^2(\real_+, Y)$,
where $Y=(\Ran V^\star)^\bot$. Since $(\Ran V^\star)^\bot=\Ker V$ is infinite dimensional, the statement follows.
$\bs$

%\subsection{The strong-star topology}\label{ssstrongstar}
\section{The strong-star topology}\label{ssstrongstar}

As we have seen in the previous section, the theory of $s$-typical properties of contractions is reduced to the study of one particular operator. The reason behind this phenomenon is that $s$-convergence does not control the adjoint, i.e.\ the function $A \mapsto A^{\star}$ is not $s$-continuous. A straightforward remedy to this problem is to refine the strong topology such that taking adjoint becomes a continuous operation. This naturally leads to the investigation of the strong-star topology.

As one may expect, the structure of an $s^{\star}$-typical contraction is more complicated than the structure of an $s$-typical contraction. We show that the theory of $s^{\star}$-typical properties of contractions can be reduced to the theories of typical properties of unitary and positive self-adjoint operators in the better understood strong topology.

\begin{thm}\label{reduc} There exist an $s^{\star}$-co-meager $s^{\star}$-$G_{\delta}$ set $\mc{H} \ss C(H)$ and an $s$-co-meager $s$-$G_{\delta}$ set $\mc{P} \ss P(H)$ such that the function $\Psi \colon U(H) \times \mc{P} \rar \mc{H}$,  $\Psi(U,P) = U \cdot P$ is a homeomorphism, where $U(H)$ and $\mc{P}$ are endowed with the strong topology and $\mc{H}$ is endowed with the strong-star topology.

Moreover, if $(\psi_{0}, \psi_{1}) \colon \mc{H} \rar U(H) \times \mc{P}$ denotes the inverse of $\Psi$, then for every $A \in \mc{H}$ and $U \in U(H)$ we have $UAU^{-1} \in \mc{H}$ and $\psi_{i}(UAU^{-1}) = U\psi_{i}(A)U^{-1}$ $(i < 2)$.
\end{thm}

Note that Theorem \ref{reduc} immediately implies the following. Recall $\mc{O}(A) = \{UAU^{-1} \colon U \in U(H)\}$.

\begin{cor}\label{numeq} For every $A \in C(H)$, $\mc{O}(A)$ is an $s^{\star}$-meager subset of $C(H)$. In particular, $s^{\star}$-typical contractions are not unitarily equivalent.
\end{cor}
\prf We follow the notation of Theorem \ref{reduc}. Let $A \in C(H)$ be arbitrary. If $A \notin \mc{H}$ then by the unitary invariance of $\mc{H}$ we have $\mc{O}(A) \ss C(H) \sm \mc{H}$, i.e.\ $\mc{O}(A)$ is $s^{\star}$-meager. So we can assume $A \in \mc{H}$.
We have \begin{multline}\notag \Psi^{-1}(\mc{O}(A)) = \{(U\psi_{0}(A)U^{-1} , U\psi_{1}(A)U^{-1}) \colon U \in U(H)\} \ss \\ \{U\psi_{0}(A)U^{-1} \colon U \in U(H)\} \times \{ U\psi_{1}(A)U^{-1} \colon U \in U(H)\} = \mc{O}(\psi_{0}(A)) \times \mc{O}(\psi_{1}(A)).\end{multline}  By Proposition \ref{weakeq}, $\mc{O}(\psi_{0}(A)) \ss U(H)$ is $s$-meager. Hence $\Psi^{-1}(\mc{O}(A))$ is $s \times s$-meager in $U(H) \times \mc{P}$. By $\Psi$ being a homeomorphism, we obtain that $\mc{O}(A)$ is $s^{\star}$-meager in $\mc{H}$ and so in $C(H)$, as well.
The further corollary that unitarily equivalent pairs are $s^{\star} \times s^{\star}$-meager in $C(H) \times C(H)$ follows as in the proof of Proposition \ref{weakeq}.$\bs$

\medskip

Similarly to our approach in the previous section, for the proof of Theorem \ref{reduc} we need to describe the mapping properties of  $s^{\star}$-typical contractions. These investigations will also help us to determine the $s^{\star}$-typical spectral properties of contractions.

%\subsubsection{Mapping properties of $s^{\star}$-typical contractions}
\subsection{Mapping properties of $s^{\star}$-typical contractions}

We prove the following.

\begin{prop}\label{surukep} The set $$\mc{T}=\{A \in C(H) \colon \forall \lambda \in \cmpl ~(\Ran (A-\lambda \cdot I) \textrm{ is dense in } H )\}$$ is $s^{\star}$-co-meager and $s^{\star}$-$G_{\delta}$ in $C(H)$.
\end{prop}

We need the following elementary property of the backward unilateral shift operator.  It follows from the fact that its adjoint, the unilateral shift operator, on $\ell^2$ has empty point spectrum (see e.g.\ \cite{R}).

\begin{lem}\label{sh1} Let $\{e_i \colon i \in \nat\}$ be an orthonormal base in $H$. Define $D \in C(H)$ by $De_{0}=0$, $De_{i+1} = e_{i}$ $(i \in \nat)$. Then for every $\lambda \in \cmpl$, $\Ran(D - \lambda \cdot \Id)$ is dense in $H$.
\end{lem}

\medskip

\textbf{Proof of Proposition \ref{surukep}.} First we show that $\mc{T}$ is $s^{\star}$-dense in $C(H)$. Let $U \ss C(H)$ be a non-empty $s^{\star}$-open set. Then there exist $\{x_{i} \colon i < n\} \ss S_{H}$, $\ve > 0$ and $A \in C(H)$ such that for every $B \in C(H)$, $\|Bx_{i} - Ax_{i}\| \leq \ve$ $(i < n)$ and $\|B^{\star}x_{i} - A^{\star}x_{i}\| \leq 2\ve$ $(i < n)$ imply $B \in U$. It is enough to find a $B \in U$ such that for every $\lambda \in \cmpl$, $\Ran (B- \lambda \cdot \Id)$ is dense in $H$.

Set $V=\span{x_i, Ax_i,
A^{\star}x_i \colon i<n}$.
Let $Q \colon V \rar V$ be defined by $Q = (1-\ve) \cdot \Pr_{V} A|_{V}$.
Let $\dim V = m$ and let $\{x_i(0) \colon i < m\}$ be an orthonormal base in $V$. Let $\{x_i(j) \colon i < m, j \in \nat\sm\{0\}\}$ be an orthonormal base in $V^{\bot}$. Define $T \colon V^{\bot} \rar H$ by $Tx_{i}(j)=x_{i}(j-1)$ $(i < m, j \in \nat \sm \{0\})$. Set $B = Q \oplus \ve \cdot T$; we show that $B$ fulfills the requirements.

We have $\|B\| \leq \|B|_{V}\|+ \|B|_{V^{\bot}}\| = \|Q\|+\ve
\cdot \|T\| \leq 1-\ve + \ve = 1$, so $B \in C(H)$. For every
$\lambda \in \cmpl$, $$\Ran (B - \lambda \cdot \Id) \supseteq \Ran
(B|_{V^{\bot}} - \lambda \cdot \Id |_{V^{\bot}}) \supseteq \Ran
\left(T - \frac{\lambda}{\ve} \cdot \Id |_{V^{\bot}}\right).$$
Since $T$ is an $m$-fold sum of the operator $D$ of Lemma
\ref{sh1}, $\Ran \left(T - \frac{\lambda}{\ve} \cdot \Id
|_{V^{\bot}}\right)$ is dense in $H$, as required.

For every $i < n$, $$\|Bx_{i} - Ax_{i}\| = \|Qx_{i} - Ax_{i}\| = \|Qx_{i} - \TS \Pr_{V} A|_{V}x_{i}\| \leq \ve,$$ \begin{multline}\notag \|B^{\star}x_{i} - A^{\star}x_{i}\| = \|Q^{\star}x_{i} + \ve \cdot T^{\star}x_{i} - A^{\star}x_{i}\| \leq \\ \ve \cdot \|T^{\star}x_{i}\| + \|Q^{\star}x_{i} - A^{\star}x_{i}\|  \leq \ve + \|Q^{\star}x_{i} - \TS \Pr_{V} A^{\star}|_{V}x_{i}\| \leq 2 \ve,\end{multline} i.e.\ $B \in U$, as required.

It remains to show that $\mc{T}$ is $s^{\star}$-$G_{\delta}$ in
$C(H)$. To this end, for every $y \in S_{H}$,
$\delta > 0$
and $L \geq 0$ set
\Keq
\label{CSI}
R(y,\delta,L) = \{A \in C(H) \colon \exists \lambda
\in \cmpl~(|\lambda| \leq L,~\dist(\Ran (A-\lambda \cdot \Id),y) \geq
\delta\}.
\Zeq
We show that $R(y,\delta,1)$ $(y \in S_{H}, \delta > 0)$ are $s$-closed. Let $\{A_n \colon n \in \nat\} \ss R(y,\delta, 1)$ and $A \in C(H)$ be such that $s$-$\lim_{n \in \nat} A_n = A$;  we show that $A \in R(y,\delta, 1).$
By $A_n \in R(y,\delta, 1)$ $(n \in \nat)$, we have $\lambda_n \in \cmpl$, $|\lambda_n| \leq 1$ $(n \in \nat)$ such that $\dist(\Ran(A_n-\lambda_n \cdot \Id),y) \geq \delta$.

By passing to a subsequence, we can assume $\{\lambda_n \colon n \in \nat\}$ is convergent, say
$\lim_{n \in \nat}\lambda_n = \lambda$; then $|\lambda|\leq 1$. It's enough to prove $\dist(\Ran(A-\lambda \cdot \Id),y) \geq \delta$. Suppose this is not the case, i.e.\ there is an $x \in H$ such that
$\|Ax-\lambda x -y\| < \delta$. Since $\lim_{n \in \nat} A_n x = A x$ and $\lim_{n \in \nat} \lambda_n =\lambda$, for an $n$ sufficiently large we have $\|A_n x -\lambda_n x -y \| < \delta$, which contradicts
$\dist(\Ran(A_n -\lambda_n \cdot \Id),y) \geq \delta$.

Thus $R(y,\delta,1)$ $(y \in S_{H}, \delta > 0)$ are $s$-closed, and so  $s^{\star}$-closed, as
well. Let $Y \ss S_{H}$ be a dense countable set. Since $\mc{T} = C(H)
\sm \TS \bigcup\{R(y,2^{-n},1) \colon y \in Y, n \in \nat\},$ the
statement follows.
$\bs$

\medskip

For technical reasons, we state the following corollary of Proposition \ref{surukep}.

\begin{cor}\label{sk_cor} The set $$\mc{E}=\{A \in C(H) \colon \Ran A \textrm{ is dense in }H \}$$  is $s^{\star}$-co-meager and $s^{\star}$-$G_{\delta}$ in $C(H)$.
\end{cor}
\prf By Proposition \ref{surukep}, $\mc{E}$ is $s^{\star}$-co-meager in $C(H)$. With the notation of (\ref{CSI}),  $\mc{E} = C(H) \sm \TS \bigcup\{R(y,2^{-n},0) \colon y \in Y, n \in \nat\},$ so the statement follows.$\bs$

\medskip

We also need a similar result for positive self-adjoint operators.

\begin{prop}\label{surukep_pos} The set $$\mc{P} = \{P \in P(H) \colon \Ran P \textrm{ is dense in } H \}$$ is $s$-co-meager and $s$-$G_{\delta}$ in $P(H)$.
\end{prop}
\prf First we show that $\mc{P}$ is $s$-dense in $P(H)$. To this end, let $U \ss P(H)$ be a non-empty $s$-open set. Then there exist $\{x_{i} \colon i < n\} \ss S_{H}$, $\ve > 0$ and $A \in P(H)$ such that for every $B \in P(H)$, $\|Bx_{i} - Ax_{i}\| \leq \ve$ $(i < n)$ implies $B \in U$. We find a $B \in U$ such that $\Ran B$ is dense in $H$.

Set $V = \span{x_{i}, Ax_{i} \colon i < n }$. Let $Q \colon V \rar V$ be an invertible contractive positive self-adjoint operator such that $\|Q - \Pr_{V} A|_{V}\| \leq \ve$; such a $Q$ exists, e.g.\ any $Q = (1-\delta)\cdot\Pr_{V} A|_{V} + \delta' \cdot \Id|_{V}$ fulfills the requirements for suitable $0 < \delta,  \delta'  \leq \ve/2$. Set $B = Q \oplus \Id|_{V^{\bot}}$. Then $B$ is a positive self-adjoint operator, $\|B|_{V}\| \leq 1$, $\|B|_{V^{\bot}}\| = 1$ and $B$ is invertible, hence $B \in P(H)$ and $\Ran B$ is dense in $H$. For every $i < n$, $$\|Bx_{i} - Ax_{i}\| = \|Qx_{i} - Ax_{i}\| = \|Qx_{i} - \TS \Pr_{V} A|_{V}x_{i}\| \leq \ve,$$ i.e.\ $B \in U$, as required.

It remains to show that $\mc{P}$ is $s$-$G_{\delta}$. Observe that for every $y \in S_{H}$ and $\delta > 0$, the set $R(y,\delta) = \{A \in P(H) \colon \dist(\Ran A,y) \geq \delta\}$ is $s$-closed. Let $Y \ss S_{H}$ be a dense countable set, then $$\mc{P} = P(H) \sm \TS \bigcup\{R(y,2^{-n}) \colon y \in Y, n \in \nat\},$$ which completes the proof.$\bs$

%\subsubsection{Proof of Theorem \ref{reduc}}
\subsection{Proof of Theorem \ref{reduc}}

With the notation of Corollary \ref{sk_cor}, set $$\mc{H} = \mc{E} \cap \{A^{\star} \colon A \in \mc{E}\} = \{A \in C(H) \colon \Ran A \textrm{ and } \Ran A^{\star} \textrm{ are dense in }H\}.$$ By Corollary \ref{sk_cor}, $\mc{H}$ is an $s^{\star}$-co-meager $s^{\star}$-$G_{\delta}$ set in $C(H)$. Define $\psi_{1} \colon \mc{H} \rar P(H)$ by $\psi_{1}(A) = (A^{\star}A)^{1/2}$ $(A \in \mc{H})$.
For every $A \in C(H)$, $(A^{\star}A)^{1/2}$ is a positive self-adjoint contraction so the definition makes sense. Moreover, for $A \in \mc{H}$ we have $\Ran (A^{\star}A)^{1/2} \supseteq \Ran (A^{\star}A)$ is dense in $H$, so with the notation of Proposition \ref{surukep_pos}, $\psi_{1}$ maps into $\mc{P}$. Note that the mapping $A \mapsto A^{\star}A$ is $s^{\star}$-continuous. We need that $\psi_{1}$ is $s^{\star}$-continuous. Since we couldn't find a reference, we outline a proof of the following.

\begin{lem}\label{sqrt} The function $\cdot^{1/2} \colon P(H) \rar P(H)$, $A \mapsto A^{1/2}$ $(A \in P(H))$ is $s$-continuous.
\end{lem}
\prf Let $A \in P(H)$ be arbitrary; note that $\sigma(A) \ss [0,1]$. By \cite[Theorem XI.6.1 p.\ 313] {Y}, \cite[Theorem XI.4.1 p.\ 307]{Y}, \cite[Proposition XI.5.2 p.\ 310]{Y} and \cite[Theorem XI.8.1 p.\ 319]{Y}, there is a $s$-left-continuous function $F \colon [0,1] \rar P(H)$ such that $F$ is projection-valued, $F(0) = 0$, $F(1) = \Id$, $F(\vartheta)F(\vartheta') = F(\min \{ \vartheta,\vartheta'\})$ $(\vartheta, \vartheta' \in [0,1])$, and for every continuous function $f \colon [0,1] \rar \real$, $\int_{0}^{1}f(\vartheta) d F(\vartheta)$ exists as the $s$-limit of Riemann-Stieltjes sums and defines a bounded linear operator, denoted by $f(A)$, satisfying $\|f(A)\| \leq \|f\|_{\infty}$. Moreover, by \cite[Example XI.12.1 p.\ 338]{Y} and \cite[Theorem XI.12.3 p.\ 343]{Y}, this defines a functional calculus compatible with power series expansion, in particular $A^{1/2}A^{1/2} = A$.

Let $\{A_{n} \colon n \in \nat\} \ss P(H)$ satisfy $s$-$\lim _{n \in \nat} A_{n} = A$. Let $x \in S_{H}$ and $\ve > 0$ be arbitrary. Let $p \colon [0,1] \rar \real$ be a polynomial satisfying $\max_{x \in [0,1]} |p(x)-x^{1/2}| < \ve$. Since $s$-$\lim _{n \in \nat} A^{k}_{n} = A^{k}$ $(k \in \nat)$, we have $\lim_{n \in \nat} p(A_{n}) x = p(A)x$. As we observed above, $\|A^{1/2} - p(A)\|, \|A^{1/2}_{n} - p(A_{n})\| < \ve$ $(n \in \nat)$. This implies $\limsup_{n \in \nat} \|A^{1/2}x - A_{n}^{1/2} x \| < 2 \ve$; since $\ve$ was arbitrary, the statement follows.$\bs$

\begin{cor}\label{sqrtcor} The function $\psi_{1} \colon \mc{H} \rar \mc{P}$ is $s^{\star}$-continuous.
\end{cor}
\prf Since the strong and strong-star topologies coincide on $P(H)$, the statement follows.$\bs$

\medskip

Let $A \in \mc{H}$ be arbitrary. As we observed above, $\psi_{1}(A)$ is a positive self-adjoint operator with dense range. Hence $\psi_{1}(A)^{-1}$ is a closed densely defined positive self-adjoint operator. Consider the densely defined operator $A \cdot \psi_{1}(A)^{-1}$. It has a dense range, and for every $x \in \Ran \psi_{1}(A)$ we have \begin{multline}\notag \sp{A \cdot \psi_{1}(A)^{-1} x,A \cdot \psi_{1}(A)^{-1}x} = \sp{A^{\star}A \cdot \psi_{1}(A)^{-1}x,  \psi_{1}(A)^{-1}x} = \\ \sp{\psi_{1}(A)^{-1} \cdot A^{\star} A \cdot \psi_{1}(A)^{-1}x, x} = \sp{x,x}.\end{multline} Hence  $A \cdot \psi_{1}(A)^{-1}$ is a densely defined isometry with dense range, i.e. it is closable and its closure is an unitary operator. Let $\psi_{0}(A) \in U(H)$ denote the closure of $A \cdot \psi_{1}(A)^{-1}$. Note that by definition, $A=\psi_{0}(A) \cdot \psi_{1}(A)$.

\begin{cor}\label{csakugy} Every $A \in \mc{H}$ can be written as $A=\psi_{0}(A) \cdot \psi_{1}(A)$ where $\psi_{0}(A) \in U(H)$ and $\psi_{1}(A) \in \mc{P}$. This decomposition is unique and unitary invariant, i.e.\ for every $A \in \mc{H}$ and $U \in U(H)$ we have $\psi_{i}(UAU^{-1}) = U\psi_{i}(A)U^{-1}$ $(i < 2)$.
\end{cor}
\prf  To see uniqueness, let $U,V \in U(H)$ and $P,Q \in \mc{P}$ be arbitrary. Then $U\cdot P= V\cdot Q$ implies $P^2 = P^{\star}U^{\star} UP = Q^{\star}V^{\star} VQ = Q^2$, so $P=Q$ by the uniqueness of positive square root. By $P,Q$ having dense range, we get  $U=V$.

If $A \in \mc{H}$ and $U \in U(H)$ then $UAU^{-1} \in \mc{H}$ and
$UAU^{-1} = U\psi_{0}(A)U^{-1} \cdot U\psi_{1}(A)U^{-1}$. Hence by the
uniqueness of this decomposition, we have $\psi_{0}(UAU^{-1})=
U\psi_{0}(A)U^{-1}$ and $\psi_{1}(UAU^{-1}) = U\psi_{1}(A)U^{-1}$.
$\bs$

\medskip

Recall $\Psi \colon U(H) \times \mc{P} \rar \mc{H},~ \Psi(U,P) = U \cdot P$, which is an $s^{\star}$-continuous function. By Corollary \ref{csakugy}, $\Psi$ is injective. It is obvious that $(\psi_{0},\psi_{1})$ is the inverse of $\Psi$, so the following completes the proof of Theorem \ref{reduc}.

\begin{lem}\label{psi1} The function $\psi_{0}$ is $s^{\star}$-continuous.
\end{lem}
\prf Let $\{A_{n} \colon n\in \nat\} \ss \mc{H}$ and $A \in \mc{H}$ be
such that $s^{\star}$-$\lim_{n \in \nat} A_{n} = A$. For every $y \in
\Ran \psi_{1}(A)$, say $y = \psi_{1}(A)x$, we have
\begin{multline}\notag\|\psi_{0}(A_{n})y - \psi_{0}(A)y\| =
  \|\psi_{0}(A_{n})\psi_{1}(A)x - \psi_{0}(A)\psi_{1}(A)x\| \leq \\
  \|\psi_{0}(A_{n})\psi_{1}(A)x - \psi_{0}(A_{n})\psi_{1}(A_{n})x\| +
  \|\psi_{0}(A_{n})\psi_{1}(A_{n})x - \psi_{0}(A)\psi_{1}(A)x\| \leq
  \\ \|\psi_{0}(A_{n})\| \cdot \|\psi_{1}(A)x - \psi_{1}(A_{n})x\| +
  \|A_{n}x - Ax\|. \end{multline} Since $\psi_{1}$ is
$s^{\star}$-continuous, the statement follows from $\Ran \psi_{1}(A)
\ss H$ being dense.$\bs$

%\subsubsection{Spectral properties of $s^{\star}$-typical contractions}
\subsection{Spectral properties of $s^{\star}$-typical contractions}

Our main result in this section is the following.

\begin{prop}\label{kisp} An  $s^{\star}$-typical contraction $A$ satisfies $C_{\sigma}(A) = \{\lambda \in \cmpl \colon |\lambda| \leq 1\}$.
\end{prop}

We start with a lemma.

\begin{lem}\label{sk2} The set $$\mc{S} = \{A \in C(H) \colon \forall \lambda \in \cmpl, |\lambda| \leq 1 ~(\lambda \in \sigma(A))\}$$ is $s^{\star}$-co-meager and $s^{\star}$-$G_{\delta}$ in $C(H)$.
\end{lem}
\prf First we show that $\mc{S}$ is $s^{\star}$-dense in $C(H)$. Let $U \ss C(H)$ be a non-empty $s^{\star}$-open set. Then there exist $\{x_{i} \colon i < n\} \ss S_{H}$, $\ve > 0$ and $A \in C(H)$ such that for every $B \in C(H)$, $\|Bx_{i} - Ax_{i}\| < \ve$ $(i < n)$ and $\|B^{\star}x_{i} - A^{\star}x_{i}\| < \ve$ $(i < n)$ imply $B \in U$. It is enough to find a $B \in U$ such that for every $\lambda \in \cmpl$ with $|\lambda| \leq 1$, $\lambda \in \sigma(B)$.

Set $V=\span{x_i,Ax_i,A^{\star}x_i \colon i<n}$, and let
$Q \colon V \rar V$ be defined by $Q = \Pr_{V} A|_{V}$. Let $T \in
C(V^{\bot})$ be arbitrary satisfying $\lambda \in \sigma(T)$ $(\lambda
\in \cmpl, |\lambda| \leq 1)$, say the backward unilateral shift operator,
and set $B = Q \oplus T$. We show that $B$ fulfills the requirements.

Since $\sigma(T) \ss \sigma(B)$, we have $\lambda \in \sigma(B)$ $(\lambda \in \cmpl, |\lambda| \leq 1)$.
For every $i < n$, $$\TS \|Bx_{i} - Ax_{i}\| = \|Qx_{i} - Ax_{i}\| = \| \Pr_{V} A|_{V}x_{i} - Ax_{i}\| =0 < \ve,$$ $$\TS \|B^{\star}x_{i} - A^{\star}x_{i}\| = \|Q^{\star}x_{i} - A^{\star}x_{i}\| = \| \Pr_{V} A^{\star}|_{V}x_{i} - A^{\star}x_{i}\| =0 < \ve,$$ i.e.\ $B \in U$, as required.

It remains to show that $\mc{S}$ is $s^{\star}$-$G_{\delta}$ in $C(H)$. For every $\delta > 0$, set $$S(\delta)= \{A \in C(H) \colon \exists \lambda \in \cmpl ~(|\lambda | \leq 1 \textrm{ and }\|(A - \lambda \cdot I)x\|, \|(A^{\star} - \overline \lambda \cdot I)x\| \geq \delta ~(x \in S_{H})) \}.$$
It is clear that $S(\delta)$ $(\delta > 0)$ are $s^{\star}$-closed. So the proof will be complete if we show $\mc{S} = C(H) \sm \bigcup_{n \in \nat}
S(2^{-n})$.

Let $A \in C(H)$ be arbitrary. If $\lambda \notin \sigma(A)$ for a $\lambda \in \cmpl$, $|\lambda| \leq 1$ then $(A - \lambda \cdot I)^{-1}$ and $(A^{\star} - \overline \lambda \cdot I)^{-1}$ exist so $A \in S(\delta)$ for some sufficiently small $\delta > 0$. This proves $C(H) \sm \mc{S} \ss \bigcup_{n \in \nat}
S(2^{-n})$.

Similarly, if $A \in S(\delta)$ for some $\delta > 0$, say $\lambda \in \cmpl$, $|\lambda| \leq 1$ witnesses this, then $\lambda \notin
P_{\sigma}(A) \cup C_{\sigma}(A)$ and $\overline \lambda \notin
P_{\sigma}(A^{\star})$.  By $\Ker (A^{\star} - \overline \lambda \cdot I) = \Ran (A - \lambda \cdot I) ^{\bot}$,      $\overline \lambda \notin
P_{\sigma}(A^{\star})$ implies $\lambda \notin R_{\sigma}(A)$. So we obtained $\lambda \notin \sigma(A)$, i.e.\ $A \notin \mc{S}$. This proves $\bigcup_{n \in \nat}
S(2^{-n}) \ss C(H) \sm \mc{S}$ and finishes the proof.$\bs$

\medskip

\textbf{Proof of Proposition \ref{kisp}.} The map $A \mapsto A^{\star}$ is a $s^{\star}$-homeomorphism of $C(H)$. So by Proposition \ref{surukep}, a typical contraction $A$ satisfies that for every $\lambda \in \cmpl$, $\Ran (A-\lambda \cdot I)$ and  $\Ran (A^{\star}-\lambda \cdot I)$ are dense in $H$. Since $\Ker (A-\lambda \cdot I) = \Ran (A^{\star}-\overline \lambda \cdot I)^{\bot}$, we get $\Ker (A-\lambda \cdot I) = \es$ $(\lambda \in \cmpl)$. Hence $P_{\sigma}(A) = P_{\sigma}(A^{\star}) = \es$ and $R_{\sigma}(A) = R_{\sigma}(A^{\star}) = \es$. So $C_{\sigma}(A) = \{\lambda \in \cmpl \colon |\lambda| \leq 1\}$ is an immediate corollary of Lemma \ref{sk2}.$\bs$

%\subsection{The norm topology}\label{ssnorm}
\section{The norm topology}\label{ssnorm}

In the norm topology it is not possible to give a simple description of the spectral properties of typical operators. An intuitive explanation may be that the norm topology is non-separable, hence several different properties can coexist on non-meager sets.

In this section, every topological notion refers to the norm topology. We prove the following.

\begin{thm}\label{mormcor} Let $\lambda \in \cmpl$ be arbitrary. Then  the following sets of operators have non-empty interior.
\ben
\item \label{mormcor_1} $\{A \in B(H) \colon \lambda \notin \sigma(A)\}$;
\item \label{mormcor_2} $\{A \in B(H)\colon  \lambda \in R_{\sigma}(A)\}$;
\item \label{mormcor_3} $\{A \in B(H) \colon  \lambda \in P_{\sigma}(A)\}$.
\een
In particular, the following sets have non-empty interior.
\ben \setcounter{enumi}{3}
\item \label{mormcor_6} $\{A \in B(H) \colon \Ran(A-\lambda \cdot I) =H\}$;
\item \label{mormcor_4} $\{A \in B(H) \colon \Ran(A-\lambda \cdot I) \textrm{ is not dense in $H$}\}$.
\een
On the other hand, the following sets of operators are nowhere dense.
\ben
\setcounter{enumi}{5}
\item \label{mormcor_7} $\{A \in B(H) \colon \lambda \in C_{\sigma}(A)\}$;
\item \label{mormcor_8} $\{A \in B(H) \colon \Ran(A-\lambda \cdot I) \textrm{ is dense in $H$ but not equal to $H$}\}$.
\een
\end{thm}

We need some terminology in advance.

\begin{defin}\rm \label{STa} Let $A \in B(H)$ and $\lambda \in \cmpl$ be arbitrary. We say $\lambda \in \sigma(A)$ is  \emph{stable} if there is an $\ve > 0$ such that for every $D \in B(H)$ with $\|D\| < \ve$ we have $\lambda \in \sigma(A+D)$. Similarly, we say $\lambda \in P_{\sigma}(A)$ is  \emph{stable} if there is an $\ve > 0$ such that for every $D \in B(H)$ with $\|D\| < \ve$ we have $\lambda \in P_{\sigma}(A+D)$.
\end{defin}

With this terminology, \cite[Theorem 2 p.\ 912]{feldman/kadison:1954} can be reformulated as follows.

\begin{prop}\label{sstabil} Let $A \in B(H)$ and $\lambda \in \cmpl$ be arbitrary.  Then the following are equivalent.
\ben
\item\label{sstabil_2} $\lambda \in \sigma(A)$ is  stable;
\item\label{sstabil_3} $\dim \Ker(A-\lambda \cdot I) \neq \dim \Ran(A-\lambda \cdot I)^{\bot}$ and  there is an $\ve > 0$ such that for every $x \in \Ker(A-\lambda \cdot I)^{\bot}$ we have $\|(A-\lambda \cdot I)x\| \geq \ve \cdot \|x\|$.
\een
\end{prop}

We prove a similar result for the point spectrum.

\begin{prop}\label{pstabil} Let $A \in B(H)$ and $\lambda \in \cmpl$ be arbitrary.  Then the following are equivalent.
\ben
\item\label{pstabil_1} $\lambda \in P_{\sigma}(A)$ is  stable;
\item\label{pstabil_3} $\dim \Ran(A-\lambda \cdot I)^{\bot} < \dim \Ker(A-\lambda \cdot I)$ and  there is a $\ve > 0$ such that for every $x \in \Ker(A-\lambda \cdot I)^{\bot}$ we have $\|(A-\lambda \cdot I)x\| \geq \ve \cdot \|x\|$.
\een
\end{prop}
\prf By linearity, we can assume $\lambda = 0$. Suppose first $0 \in P_{\sigma}(A)$ is  stable. Then $0 \in \sigma(A)$ is stable so by Proposition \ref{sstabil}, $\dim \Ker(A) \neq \dim \Ran(A)^{\bot}$ and  there is a $\ve > 0$ such that for every $x \in \Ker(A)^{\bot}$ we have $\|Ax\| \geq \ve \cdot \|x\|$. If  $\dim \Ker(A) < \dim \Ran(A)^{\bot}$ then for every $\ve > 0$ there is a $D \in B(H)$ such that $\|D\| < \ve$, $D|_{\Ker(A)^{\bot}} =0$ and for every $x \in \Ker(A) \sm \{0\}$, $Dx \in  \Ran(A)^{\bot} \sm \{0\}$. Then for every $x \in H\sm \{0\}$ we have $(A+D)x \neq 0$, hence $0 \notin P_{\sigma}(A+D)$. This contradicts the assumption that $0 \in P_{\sigma}(A)$ is stable. Thus $\dim \Ran(A)^{\bot} < \dim \Ker(A)$, as required.

To see the converse, suppose the conditions of statement \ref{pstabil_3} hold. Then we have $\dim \Ker(A) > 0$ hence $0 \in P_{\sigma}(A)$. We also have that $A \colon \Ker(A)^{\bot} \rar \Ran(A)$ is a bounded invertible operator, in particular $\Ran(A)$ is a closed co-finite dimensional subspace of $H$. Let $k = 1+\dim \Ran(A)^{\bot}$.

Let $D \in B(H)$ be arbitrary with $\|D\| < \ve/2k^2$. For every $x \in \Ker(A)$ consider the following inductive definition of a sequence $\{x_{n} \colon n \in \nat\} \ss H$. Set $x_{0} = x$. Let $n \in \nat$ be arbitrary and suppose that $x_{n} \in H$ is defined. Write $Dx_{n} = u_{n} + v_{n}$ where $u_{n} \in \Ran(A)$ and $v_{n} \in \Ran(A)^{\bot}$. Set $x_{n+1} = A^{-1} u_{n}$. This completes the inductive step of the definition of $\{x_{n} \colon n \in \nat\}$.

We set $\xi(x) = \sum_{n \in \nat} (-1)^{n}x_{n}$, $\rho(x) = \sum_{n \in \nat} (-1)^{n}v_{n}$. Note that $\|x_{n+1}\| \leq \|u_{n}\|/\ve \leq \|Dx_{n}\|/\ve < \|x_{n}\|/2k^2$, hence $\|x_{n}\| < \|x\|/(2k^{2})^{n}$ $(n \in \nat)$ so the definitions make sense. Moreover, the functions $\xi$ and $\rho$ are linear, $\rho(x) \in \Ran(A)^{\bot}$ and $\|\xi(x)-x\| < 1/k$ $(x \in B_{\Ker(A)})$. So by Lemma \ref{ortho}, if $\{x(i) \colon i  < k\}$ is an orthonormal system in $\Ker(A)$ then $\{\xi(x(i)) \colon i  < k\}$ are linearly independent.

Observe that \begin{multline}\notag  (A+D)\xi(x) = (A+D)  \sum_{n \in \nat} (-1)^{n}x_{n} =  \sum_{n \in \nat} (-1)^{n}
(Ax_{n} + Dx_{n}) = \\ Ax_{0} +  \sum_{n \in \nat \sm \{0\}} (-1)^{n}Ax_{n} + \sum_{n \in \nat} (-1)^{n}(u_{n}+v_{n}) = \sum_{n \in \nat} (-1)^{n}v_{n} = \rho(x). \end{multline}
By $\dim \Ker(A) \geq k$, there is an orthonormal system  $\{x(i) \colon i < k\}$ in $\Ker(A)$. Since $\dim \Ran(A)^{\bot} < k$, there is an $x \in \span{x(i) \colon i < k} \sm \{0\}$ such that $\rho(x)=0$. Then $(A+D)\xi(x)=0$ shows $0 \in P_{\sigma}(A+D)$, so the proof is complete.$\bs$

\medskip

\textbf{Proof of Theorem \ref{mormcor}.} By linearity, it is enough
to consider the $\lambda = 0$ case. The set of invertible operators
shows statements \ref{mormcor_1} and \ref{mormcor_6}. By Proposition
\ref{pstabil}, a neighborhood of the
backward unilateral shift operator shows statement \ref{mormcor_3}.

The set of statement \ref{mormcor_8} contains the set of statement
\ref{mormcor_7} so it is sufficient to prove statement
\ref{mormcor_8}. Let $A \in B(H)$ be arbitrary satisfying $\Ran(A)$ is
dense in $H$ but not equal to $H$. Then $A$ cannot be
invertible. Moreover, note that condition $\|Ax\| \geq \ve \cdot
\|x\|$ $(x \in \Ker(A)^{\bot})$ in Proposition
\ref{sstabil}.\ref{sstabil_3} would imply $\Ran(A)$ is closed.
Since $\Ran(A)$ cannot be closed, $A$ cannot satisfy Proposition \ref{sstabil}.\ref{sstabil_3}.
So the set of statement \ref{mormcor_8} is contained in the boundary of the set of invertible operators hence it is nowhere dense.

Finally notice that a neighborhood of the unilateral shift operator $D^{\star}$ is contained in $\{A \in B(H)\colon  0 \in \sigma(A) \sm P_{\sigma}(A)\}$, since the stability of $0 \in \sigma(D^{\star})$ is implied by Proposition \ref{sstabil} while the stability of $0 \notin P_{\sigma}(D^{\star})$ follows from $D^{\star}$ being an isometry.  So by statement \ref{mormcor_7}, statements \ref{mormcor_2} and \ref{mormcor_4} follow.$\bs$

\medskip

Observe that the backward unilateral shift operator $D$ of Lemma \ref{sh1} satisfies the conditions of statement \ref{pstabil_3} of Proposition \ref{pstabil} for every $\lambda \in \cmpl$ with $|\lambda| < 1$, hence every such $\lambda \in P_{\sigma}(D)$ is stable. So by taking direct sums of scaled and translated copies of $D$ we can obtain arbitrarily complicated stable spectra. We also obtain the following. Here an operator $T$ is called \emph{embeddable}, if it can be embedded into a strongly continuous semigroup, i.e.~if $T=T(1)$ holds for some strongly continuous semigroup $(T(t))_{t\geq 0}$.

\begin{cor}\label{normein} The set of embeddable operators and the set of non-embeddable operators have non-empty interior.
\end{cor}
\prf By \cite[Theorem 2.1 p.\ 452]{E2}, each operator in a small
neighborhood of $\Id$ is embeddable. To show a non-empty  open set of
non-embeddable operators, we aim to use \cite[Theorem 3.1 p.\
454]{E2}. Recall the
backward unilateral shift operator $D$ of Lemma \ref{sh1}. By Proposition \ref{pstabil}, each operator $A$ in a small neighborhood of $D$ satisfies $0 \in P_{\sigma}(A)$. So we get $A$ is non-embeddable if we show $\dim \Ker (A) \leq 1$.

Note that if $\|A - D\| < 1/\sqrt{2}$ then for every $x,y \in \Ker (A) \cap S_{H}$ we have $\|Dx\|, \|Dy\| < 1/\sqrt{2}$. This implies $\|x - D^{\star}Dx\|, \|y - D^{\star}Dy\| >  1/\sqrt{2}$, and since $D^{\star}D$ is an orthogonal projection, we have
\begin{multline}\notag |\sp{x,y}| = |\sp{x - D^{\star}Dx, y - D^{\star}Dy} + \sp{D^{\star}Dx,D^{\star}Dy}| \geq \\
|\sp{x - D^{\star}Dx, y - D^{\star}Dy}| - |\sp{D^{\star}Dx,D^{\star}Dy}| > 0.
\end{multline}
Thus $\dim \Ker (A) \leq 1$, which completes the proof.
$\bs$

\medskip

The results of this section provide no information about the size of the spectrum of a typical operator. We state two questions related to this in Problem \ref{PR2}. Here we only point out the following.

\begin{prop}\label{esru} The following sets are dense.
\ben
\item\label{mormcor_9} $\{A\in B(H) \colon P_\sigma(A)\neq \emptyset\}$;
\item\label{mormcor_10} $\{A\in B(H) \colon R_\sigma(A)\neq \emptyset\}$.
\een
\end{prop}
\prf Since $A \mapsto A^{\star}$ $(A \in B(H))$ is a homeomorphism of $B(H)$ and $P_\sigma(A^{\star})=R_\sigma(A)$ $(A \in B(H))$, it suffices to show \ref{mormcor_9}.

Let $B \in B(H)$ be arbitrary. Since the approximative point spectrum of $B$ is nonempty (see e.g. \cite[Proposition VII.6.7 p.\ 210]{Conway}), there exists a $\lambda\in \sigma(B)$ and a sequence $\{x_n \colon n \in \nat \} \subseteq S_H$ such that $\lim_{n\in \nat}\|Bx_n-\lambda x_n\|=0$. For every $n \in \nat$, let $P_n$ denote the orthogonal projection onto $\span{x_n}$, and set $A_n =B(\Id - P_n)+\lambda  \cdot P_n$. We have $A_{n} x_n=\lambda x_n$ and $\|A_n-B\|=\|\lambda \cdot P_n -BP_n\| = \|\lambda x_n-Bx_n\|$ $(n \in \nat)$. Hence $P_{\sigma}(A_{n}) \neq \es$ $(n \in \nat)$ and $\lim_{n \in \nat} \|A_n-B\| =  0$. This completes the proof.$\bs$

\section{An outlook to general Banach spaces}\label{sgenB}

The proofs in the previous section could convince the reader that an attempt to extend our investigations to tackle the typical properties of contractions on arbitrary Banach spaces could encounter considerable technical difficulties. In addition, such an endeavor has to face some problems of more fundamental nature. Recent developments in the theory of Banach spaces resulted in numerous spaces exhibiting surprising functional analytic properties. On the famous Banach space of S.\ A.\ Argyros and R.\ G.\ Haydon, every bounded linear operator is of the form $\lambda \cdot \Id + K$ where $\lambda \in \cmpl$ and $K$ is a compact operator (see \cite{AH}). Spaces were constructed with only trivial isometries (see e.g.\ \cite{DK} and the references therein), moreover even a renorming of a Banach space may result in an arbitrary isometry group (see \cite{FG} and the references therein). The relevance of isometries comes from the important role played by unitary operators in the theory developed in the previous section.

Observe that renorming does not change the topology of the underlying Banach space, so the five topologies we consider on operators remain unchanged, as well. Instead, renorming affects the size of various classes of operators. So in a sense, the study of typical properties of operators in various topologies is more related to the geometry of the underlying Banach space than to the topology it carries. Therefore, reasonable extensions of our investigations should be pursued in Banach spaces where the geometry of the space is of special significance.

The most obvious such spaces are the $L^{p}$ spaces $(1 \leq p \leq \infty)$. The functional analysis of these spaces is well-developed, isometric operators are characterized (see e.g.\ \cite{L}). Nevertheless, we expect that none of our main results extends to arbitrary $L^{p}$ spaces.

Another promising extension of the theory of typical behavior of operators could tackle Banach$^{\star}$-algebras. Developing our results to that generality could separate typical properties which are of operator theoretic nature from typical properties exploiting the geometry of the underlying space.

Finally let us propose some concrete problems which stem from the results of the previous section.

\begin{prob}\label{PR1} We say $A,B \in B(H)$ are \emph{similar} if there is an invertible operator $T \in B(H)$ such that $A = TBT^{-1}$.
\ben
\item Is it true that $w$-typical contractions are similar?
\item Is it true that $s^{\star}$-typical contractions are similar?
\een
\end{prob}

\begin{prob}\label{PR2} Determine the $s$-typical properties of self-adjoint and positive self-adjoint operators.
\end{prob}

\begin{prob}\label{PR2_5} Is  a $s^{\star}$-typical contraction embeddable into a strongly continuous semigroup?
\end{prob}

\begin{prob}\label{PR3} ~
\ben
\item Is the set $\{A \in B(H) \colon P_{\sigma}(A) \neq \es\}$ co-meager in the norm topology?
\item Is the set $\{A \in B(H) \colon C_{\sigma}(A) = \es\}$ co-meager in the norm topology?
\een
\end{prob}

It would be instructive to examine whether suitable analogues of our results hold for strongly continuous
semigroups instead of single operators. Let $C^c(H)$ denote the set of contractive $C_0$-semigroups. Given any metric $d$ on $C(H)$, one can endow $C^c(H)$ with the  topology generated by uniform  $d$-convergence on compact
time intervals. This topology is induced by the metric
$$ \TS
d_c(T(\cdot),S(\cdot))=\sum_{n \in \nat} 2^{-n}  \sup_{t\in [0,n]} d(T(t),S(t)).$$

\begin{prob}\label{PR0}
Under which of these topologies is $C_c(H)$ a Baire space? What are the typical properties of contractive $C_0$-semigroups in these topologies?
\end{prob}

Recall that the strong topology case was treated in Section \ref{sstrong-cont}.
Some results related to Problem \ref{PR0} for the weak topology can be
found in \cite[Section 4]{E0}, \cite{ES1} and \cite{ES2}.

\end{document}